\documentclass[english]{jnsao}
\usepackage[utf8]{inputenc}
\usepackage[nameinlink,capitalize]{cleveref}
\AddToHook{env/proposition/begin}{\crefalias{theorem}{proposition}}
\AddToHook{env/lemma/begin}{\crefalias{theorem}{lemma}}
\AddToHook{env/corollary/begin}{\crefalias{theorem}{corollary}}
\AddToHook{env/definition/begin}{\crefalias{theorem}{definition}}
\AddToHook{env/example/begin}{\crefalias{theorem}{example}}
\AddToHook{env/remark/begin}{\crefalias{theorem}{remark}}
\AddToHook{env/remark/begin}{\crefalias{theorem}{remark}}
\AddToHook{cmd/appendix/before}{\crefalias{section}{appendix}}

\usepackage{graphicx}
\usepackage{booktabs}
\usepackage{enumitem}

\manuscriptsubmitted{2025-05-14}
\manuscriptaccepted{2026-06-16}
\manuscriptvolume{6}
\manuscriptnumber{15673}
\manuscriptyear{2026}
\manuscriptdoi{10.46298/jnsao-2026-15673}

\title{A constructive approach to strengthen algebraic descriptions of function and operator classes}
\shorttitle{Algebraic strengthening of function/operator classes}

\author{Anne Rubbens\thanks{ICTEAM Institute, UCLouvain, Belgium. \email{anne.rubbens@uclouvain.be}, \email{julien.hendrickx@ucloivain.be}} \and Julien M. Hendrickx\footnotemark[1] \and Adrien Taylor\thanks{INRIA, \'Ecole normale sup\'erieure, CNRS, PSL Research University, France. \email{adrien.taylor@inria.fr}}}
\shortauthor{Rubbens, Hendrickx, Taylor}

\acknowledgements{
    A.~Rubbens is supported by a FNRS fellowship. J.~Hendrickx is supported by the SIDDARTA Concerted Research Action (ARC) of the Federation Wallonie-Bruxelles. A.~Taylor is supported by the European Union (ERC grant CASPER 101162889). The French government also partly funded this work under the management of Agence Nationale de la Recherche as part of the ``France 2030'' program, reference ANR-23-IACL-0008 (PR[AI]RIE-PSAI).  Views and opinions expressed are however those of the authors only.
}
\newcommand{\N}{\mathbb{N}}
\newcommand{\F}{\mathcal{F}}

\newcommand{\Op}{\mathcal{Q}_{\mu,L}}
\newcommand{\Opi}{\mathcal{Q}_{\mu,\beta}}
\newcommand{\R}{\mathbb{R}}
\newcommand{\pL}{p_{\text{\L{}}_{\mu,L}}}
\newcommand{\tpL}{\tilde{p}_{\text{\L{}}_{\mu,L}}}
\newcommand{\pblock}{p_{0,\mathbb{L}}}
\newcommand{\tpblock}{\tilde{p}_{0,\mathbb{L}}}
\newcommand{\pop}{q_{\mu,L}}
\newcommand{\tpop}{\tilde{q}_{\mu,L}}
\newcommand{\popi}{q_{\mu,\beta}}
\newcommand{\tpopi}{\tilde{q}_{\mu,\beta}}

\renewcommand{\leq}{\leqslant}
\renewcommand{\geq}{\geqslant}
\renewcommand{\succeq}{\succcurlyeq}
\renewcommand{\preceq}{\preccurlyeq}
\newcommand{\FQrestr}{\F_p(\bar Q)|_Q}

\begin{document}
\maketitle

\begin{abstract}
    It is well known that functions (resp. operators) satisfying a property~$p$ on a subset~$Q\subset \R^d$ cannot necessarily be extended to a function (resp. operator) satisfying~$p$ on the whole of~$\R^d$. Given $Q \subseteq \R^d$, this work considers the problem of obtaining necessary and ideally sufficient conditions to be satisfied by a function (resp. operator) on $Q$, ensuring the existence of an extension of this function (resp. operator) satisfying $p$ on $\R^d$.

    More precisely, given some property $p$, we present a refinement procedure to obtain stronger necessary conditions to be imposed on $Q$. This procedure can be applied iteratively until the stronger conditions are also sufficient. We illustrate the procedure on a few examples, including the strengthening of existing descriptions for the classes of smooth functions satisfying a \L{}ojasiewicz condition, convex blockwise smooth functions, Lipschitz monotone operators, strongly monotone cocoercive operators, and uniformly convex functions.

    In most cases, these strengthened descriptions can be represented, or relaxed, to semi-definite constraints, which can be used to formulate tractable optimization problems on functions (resp. operators) within those classes.
\end{abstract}

\section{Introduction}
This work addresses the question of extending functions and operators within certain classes. Before formally stating this problem in \cref{sec:formal_extension}, we informally describe it, considering function classes only. Informally, let $\F_p(\R^d)$ be a set of functions $f:\R^d\rightarrow \R\cup\{\infty\}$ satisfying some algebraic property ``$p$'' everywhere. For instance, we say $f$ belongs to the set of convex functions $\F_{p_{0,\infty}}(\R^d)$ if
\begin{equation*}
    f(\alpha x+(1-\alpha) y)\leq \alpha f(x)+(1-\alpha) f(y), \ \forall x,y \in \R^d, \ \forall \alpha \in [0,1]. \tag{$p_{0,\infty}$}
\end{equation*}
There are generally several ways to describe a set of functions (i.e., via different properties). In the case of convex functions, one can, for instance, define the class via subdifferentials \cite[Section 23]{rockafellar1997convex}
\begin{equation*}
    f(y)\geq  f(x) + \langle g_x,y-x\rangle \quad\forall y\in\mathbb{R}^d,\ x\in \R^d,\ g_x\in \R^d: \ g_x\in\partial f(x) .\tag{$p'_{0,\infty}$}
\end{equation*}
In this context, given a set $Q \subseteq \R^d$, we seek to answer the following question.
\begin{center}
    \noindent \textbf{What conditions must a function $f$ satisfy on~$Q$ to ensure existence of an \emph{extension} of $f$ that is identically equal to $f$ on~$Q$ and belongs to $\F_p(\R^d)$?}
\end{center}

We refer to such necessary and sufficient conditions as \emph{interpolation/extension conditions}. A necessary extension condition for $\F_p(\R^d)$ is $p$ itself, which is satisfied on~$Q$ by functions that admit an extension to~$\R^d$. However, this is not always sufficient since $p$ can be satisfied on $Q$ by functions that are not consistent with any function in $\F_p(\R^d)$, see, e.g.,~\cite[Figure 1]{taylor2017smooth} for the case of smooth convex functions.

In addition, several properties $p^{(1)},\ p^{(2)}, \ ...\ $ such that $\F_{p^{(1)}}(\R^d)=\F_{p^{(2)}}(\R^d)=...\ $ are generally not equivalent when imposed on sets $Q\subset \R^d$. We refer to such a condition as being \emph{stronger} than another if the set of functions satisfying it on~$Q$ is contained within the set of functions satisfying the other condition on Q. Given a property $p$, the strongest property $p'$ such that $\F_{p'}(\R^d)=\F_p(\R^d)$ is a sufficient condition for extension to $\mathcal{F}_p(\R^d)$, satisfied on $Q$ exactly by functions admitting an extension to $\R^d$.

We illustrate the different stages of this work via the running example of smooth convex functions.
\paragraph{Example.}  The class $\F_{{0,L}}(\R^d)$ of $L$-smooth convex functions on $\R^d$ can be described in numerous equivalent ways, see, e.g.,~\cite[Section 2.1.1]{nesterov2018lectures}. Arguably the most natural one consists in juxtaposing a convexity condition with Lipschitz continuity of the gradients, thereby obtaining a condition $p_{0,L}$ to be satisfied on $\R^d$  by $L$-smooth convex functions $f$, in the sense that $p_{0,L}$ is negative, $\forall x,y\in \R^d$
\begin{equation}\label{eq:convexity_gradient_lipschitz}\tag{$p_{0,L}$}
    \begin{aligned}
        p_{0,L}((x,f(x),\nabla f(x)),(y,f(y),\nabla f(y))) \leq 0 \Leftrightarrow \left\{\begin{array}{l}
                f(y)-f(x)+\langle \nabla f(y),x-y\rangle  \leq 0  \\[.1cm]
                \|\nabla f(x)-\nabla f(y)\|-L \|x-y\|    \leq 0.
        \end{array} \right.
    \end{aligned}
\end{equation}

Another description of $\F_{{0,L}}(\R^d)$ consists in juxtaposing convexity and another possible definition of smoothness, obtained by integration of the Lipschitz continuity of the gradients, see, e.g., \cite[Theorem 2.1.5]{nesterov2018lectures}, and imposing, $\forall x,y\in \R^d$,
\begin{equation}\label{eq:convexity_quadratic_UB}\tag{$p'_{0,L}$}
    \begin{aligned}
        {p}'_{0,L}((x,f(x),\nabla f(x)),(y,f(y),\nabla f(y)))\leq 0 \Leftrightarrow \left\{\begin{array}{l}
                f(y)- f(x)+\langle \nabla f(y),x-y\rangle \leq 0    \\[.1cm]
                f(x)-f(y)-\langle \nabla f(y),x-y\rangle-\tfrac{L}{2}\|x-y\|^2\leq 0.
        \end{array} \right.
    \end{aligned}
\end{equation}

By combining the two inequalities in \eqref{eq:convexity_quadratic_UB}, one can further show that $\F_{{0,L}}(\R^d)$ can be equivalently defined by the following single condition $p''_{0,L}$, see, e.g., \cite[Theorem 2.1.5]{nesterov2018lectures}, to be imposed on all~$x,y\in \R^d$
\begin{equation}\label{eq:true_smooth_convexity}\tag{$p''_{0,L}$}
    {{p}}''_{0,L}((x,f(x),\nabla f(x)),(y,f(y),\nabla f(y)))\leq 0 \,\Leftrightarrow \, f(x)- f(y)-\langle \nabla f(y),x-y\rangle-\tfrac{\|\nabla f(x)-\nabla f(y)\|^2}{2L}\leq 0.
\end{equation}

These three definitions are equivalent when imposed on any pair in $\R^d$ and are thus all necessary extension conditions when imposed on sets $Q\subset \R^d$, i.e., on any pair $x,y\in Q$. However, they are not equivalent anymore when imposed on such smaller sets. In particular, \eqref{eq:true_smooth_convexity} is stronger than both~\eqref{eq:convexity_quadratic_UB} and~\eqref{eq:convexity_gradient_lipschitz}, since one can find counterexamples which, on some sets $Q\subset \R^d$, satisfy~\eqref{eq:convexity_quadratic_UB} or~\eqref{eq:convexity_gradient_lipschitz}, but do not admit an extension satisfying $p$ to~$\R^d$~\cite{taylor2017smooth}. These two definitions are thus not sufficient extension conditions. On the contrary, it has been shown~\cite[Theorem 4]{taylor2017smooth} (see also \cite{azagra2017extension,daniilidis2018explicit} for alternate viewpoints) that \eqref{eq:true_smooth_convexity} is a necessary and sufficient extension property for $\F_{{0,L}}(\R^d)$, i.e., for any $Q\subset\R^d$, functions satisfying~\eqref{eq:true_smooth_convexity} on $Q$ are exactly those admitting a smooth convex extension to $\R^d$.

\paragraph{Motivation.} Extension conditions are particularly relevant for problems where functions in some class $\F_p(\R^d)$ are variables to certain optimization, (see, e.g.,~\cite{paty2020regularity,taylor2017smooth}), or regression problems, (see, e.g.,~\cite{hannah2013multivariate,lambert2004finite}). In this situation, one typically only needs or has access to evaluations of the function at a discrete set of points $Q=\{x_1,x_2,\ldots,x_N\}$, and requires this dataset to be consistent with a function in $\F_p(\R^d)$. It is then often easier to formulate conditions of existence of such functions in $\F_p(\R^d)$ by imposing conditions on $Q$ rather than on the full space~$\R^d$, which would imply dealing with infinite-dimensional variables. Among others, this kind of discretization technique became popular for obtaining worst-case performance bounds for optimization methods in a computer-aided way~\cite{drori2014performance,taylor2017exact}.  In this context, it has been shown that using conditions that guarantee extensibility/interpolation allows computing tight performance bounds. The same motivation holds for operators; see e.g., \cite{ryu2020operator}.

\subsection{Related works}\label{sec:related_works}
\paragraph{Extension/interpolation conditions.} The interest in extension conditions can be traced back to at least~1934, with~\cite{whitney1934differentiable} asking the question of whether or not a function with given properties, e.g., continuity, could be extended outside of its domain. The properties are formulated qualitatively, i.e., it is not required for the extension to satisfy the property with the same parameters as the initial function. In this context, extension conditions were obtained, e.g., for non-negative functions \cite{fefferman2017interpolation}, Lipschitz continuous functions~\cite{gruyer2009minimal} and convex functions~\cite{azagra2017extension,azagra2017whitney,azagra2019smooth,yan2012extension}.

The quest for interpolation conditions to obtain exact performance guarantees on optimization methods started in 2017 with smooth (weakly) convex functions~\cite{taylor2017exact,taylor2017smooth}. In this context, initial functions and extensions are required to satisfy properties with the same parameters. Building on this work, such conditions were then obtained, e.g., for smooth Lipschitz functions~\cite{thesis}, indicator functions of (strongly) convex sets~\cite{luner2024performance,thesis,taylor2017exact}, difference-of-convex functions and relatively smooth and convex functions~\cite{dragomir2021optimal}, and (strongly) monotone, cocoercive, Lipschitz or linear operators~\cite{bousselmi2023interpolation,ryu2020operator}, and for a variety of other classes of functions (see, e.g.,~\cite{guille2022gradient,goujaud2022optimal} and a partial list provided in the \href{https://pepit.readthedocs.io/en/latest/api/functions_and_operators.html}{PEPit documentation}~\cite{goujaud2024pepit}). A comparison of the different versions of an interpolation condition arising from these contexts is provided in \cite[Section 2]{rubbens2024constraint}. Natural consequences of knowing (or not knowing) interpolation conditions in the context of first-order optimization are discussed at length in, e.g.,~\cite{thesis,goujaud2023fundamental}.

To this date, however, obtaining extension conditions for specific function (resp.\!~operator) classes remains a challenge. In particular, there exists no principled way (i)~to strengthen a non-sufficient candidate extension condition or to prove it is sufficient, nor (ii)~to combine extension conditions associated with several properties, as showcased in~\cite[Section 2.2]{taylor2017smooth} and~\cite[Section 2.2]{ryu2020operator}.
\paragraph{Positive polynomials: sum of squares (SOS) and co-positive formulations.} Throughout this work, we propose descriptions of classes that can be expressed as \emph{convex} conditions to be imposed on a dataset (discrete version of functions/operators) and reformulate/relax these conditions to obtain \emph{tractable} forms for them. In particular, for some $k,n\in\mathbb{N}$, these convex conditions are often expressed as $n$-dimensional polynomials of degree $k$: $P(\theta)\in \R_k[\theta_1,\ldots,\theta_n]$, whose coefficients depend on the dataset, and which must be non-positive for all $\theta\geq 0$. Representing convex expressions, and in particular (co)-positive polynomials, in tractable/intractable ways has been a long-lasting subject of interest, and we quickly mention a few key techniques.

Setting $\gamma=\theta^2$, enforcing $P(\theta)\leq 0, \ \forall \theta\geq 0$ is equivalent to enforcing a $n$-dimensional polynomial of degree $2k$, $\bar{P}(\gamma)\in \R_{2k}[\gamma_1,\ldots,\gamma_n]$ to be non-positive everywhere. Interest in non-negative polynomials can be traced back to 1888 \cite{hilbert1888darstellung}, and they have been widely studied since; see, e.g., \cite{lasserre2009moments,marshall2008positive} for surveys. While the set of non-negative polynomials is non-tractable in general, it can be relaxed in several tractable ways, e.g., by requiring polynomials to be \emph{sum of squares} (SOS). That is, there exists a matrix $M\succeq0$ such that $\bar P(\gamma)=-X^T M X$, where $X$ is a vector of size $\binom{k+n}{n}$ containing all monomials of degree $\leq k$ in $\gamma$. It is well-known that for the special cases $(n=1,k\in\mathbb{N})$, quadratic polynomials and $(n=2,k=2)$, non-negative and SOS polynomials are equivalent \cite{hilbert1888darstellung}. In any case, the set of SOS polynomials provides a natural path to tractable convex relaxations of positive polynomials in general. Alternatively, when $n=1$, one can use the Markov-Lukasc theorem, see, e.g., \cite[Theorem 1.21]{szeg1939orthogonal}, to express non-negativity of polynomials as SDPs.

Finally, in some cases, enforcing $P(\theta)\leq 0, \ \forall \theta\geq 0$ can be equivalently formulated as requiring existence of a \emph{co-positive} matrix $M$ (i.e., $x^TMx\geq 0, \ \forall x\geq0$), such that $P(\theta)=-X(\theta)^T M X(\theta)$ for some vector $X(\theta)$. Co-positive programming has been widely investigated since 1952 \cite{motzkin1952copositive}, see, e.g., \cite{dur2010copositive} for a survey. While certifying the co-positivity of a matrix is generally $NP$-hard, it is semi-definite representable for small problem dimensions ($\leq 4$, see, e.g.,~\cite{maxfield1962matrix}) and can naturally be relaxed to tractable convex representations otherwise.

\subsection{Contributions}
We propose a systematic approach for searching extension conditions for given classes $\F_p(\R^d)$, starting from $p$ itself. To do so, we constructively and iteratively \emph{strengthen} $p$ via a refinement procedure. That is, we iteratively propose alternative conditions $p^{(1)}, \ p^{(2)}, \ ...$ that are equivalent to $p$ on $\R^d$ but stronger than $p$ on subsets $Q\subseteq \R^d$. Those alternatives become stronger and stronger as we proceed (i.e., the sets of functions satisfying them get smaller). This technique naturally complements the constructive approach to counterexamples proposed in~\cite{rubbens2024constraint}.

The contribution of this work is twofold. First, \cref{sec:iterative_str} presents the refinement technique, and illustrates it on the class of smooth convex functions, providing an original approach (not relying on integration arguments) to obtain~\eqref{eq:true_smooth_convexity} from~\eqref{eq:convexity_gradient_lipschitz}. Then, \cref{sec:applications} exploits the technique to obtain new conditions for function and operator classes that have sparked interest in the optimization community (see, e.g.,~\cite{abbaszadehpeivasti2022conditionslinearconvergencegradient,bolte2010characterizations,drori2018properties,karimi2016linear,nesterov2012efficiency,nesterov2015universal,nesterov2018lectures,ryu2020operator}) and for which no extension condition yet exists. These refined conditions are obtained via a single round of the refinement procedure, and are hence a priori only necessary extension conditions. However, we show on several examples how they quantitatively improve upon classical descriptions. In addition, we provide, for most cases, semi-definite programs (SDP) representations (in terms of the conditions imposed on a dataset) of those refined conditions. In particular, we strengthen existing descriptions of
\begin{itemize}
    \item smooth functions satisfying a \L{}ojasiewicz condition, see, e.g.,  \cite{bolte2007lojasiewicz,bolte2010characterizations,bolte2017error,lojasiewicz1963propriete,polyak1963gradientMF} for motivations, and \cite{abbaszadehpeivasti2022conditionslinearconvergencegradient} for recent analyses of the gradient method on this class: SDP formulation (\cref{sec:loja}),
    \item convex blockwise smooth functions, see, e.g., \cite{beck,nesterov2012efficiency} for motivations, and \cite{abbaszadehpeivasti2022convergence,kamri2023worst} for recent analyses of the cyclic or randomized blockwise coordinate descent methods on this class: SDP formulation (\cref{sec:blockwise}),
    \item Lipschitz strongly monotone operators, see, e.g., \cite{bauschke2017convex,ryu2020operator} for motivations, and \cite{gorbunov2022extragradient,gorbunov2022last,korpelevich1976extragradient} for recent analyses of the optimistic and extra gradients on this class: SDP formulation (\cref{sec:operators}),
    \item strongly monotone cocoercive operators, see, e.g., \cite{bauschke2017convex,ryu2020operator} for motivations: SDP formulation (\cref{sec:operators2}),
    \item uniformly convex functions, see, e.g., \cite{zǎlinescu1983uniformly} for an introduction and \cite{iouditski2014primal,nesterov2015universal}  for applications (\cref{sec:unif_conv}),
    \item smooth convex functions defined on convex subsets of $\R^d$, i.e., it is not required for the extension to be defined everywhere, see, e.g., \cite{drori2018properties} for a study of this case (\cref{sec:extension_subset}).
\end{itemize}
\subsection{Notation}
We consider the standard Euclidean space $\mathbb{R}^d$ with standard inner product $\langle \cdot,\cdot\rangle$ and induced norm $\|\cdot\|$. Given $N \in \mathbb{N}$, we let $\llbracket N\rrbracket=\{1, \cdots,N\}$. Given a set $Q\subseteq\R^d$, we denote by $|Q|$ its cardinality, where $|Q|=|\mathbb{N}|=\aleph_0$ for $Q$ countably infinite, $|Q|=|\R^d|=\mathfrak{c}$ for $Q$ uncountable, and $K\leq \aleph_0\leq\mathfrak{c}$, $\forall K\in \mathbb{N}$. Given a symmetric matrix $M\in \R^{d\times d}$, we say $M\succeq 0$ if it is positive semi-definite. Given a convex function $f:\R^d\longrightarrow \R$, we refer by $\partial f(x)$ to the subdifferential \cite{rockafellar1997convex} of $f$ at $x$, defined as the set of vectors $v_x\in \R^d$ satisfying $f(y)\geq f(x)+\langle v_x,y-x\rangle$. When clear from the context, we denote by $v_x$ a subgradient of $f$ at $x$ (i.e., $v_x\in\partial f(x)$).
\section{Iterative strengthening of a property}\label{sec:iterative_str}
\subsection{Formal extension property}  \label{sec:formal_extension}
This work considers properties that can be translated into one or several inequalities involving a finite number of points. Such properties can always be reduced to a single inequality, taking the supremum/maximum of the algebraic expressions among all points, see, e.g.,~\eqref{eq:convexity_gradient_lipschitz}, \eqref{eq:convexity_quadratic_UB} or \eqref{eq:true_smooth_convexity}.

Formally, let $Q\subseteq \R^d$ and let $F:Q\rightrightarrows \R^D$ a multi-valued mapping, and consider properties $p:(\R^d\times\R^D)^m \to \R$. We say $F$ satisfies $p$ at some tuple in $Q^m$ if $p((x_1,F_{x_1}),\ldots ,(x_m,F_{x_m}))\leq 0$, for all $F_{x_i}\in F(x_i)$, $i\in \llbracket m\rrbracket$.
\paragraph{Mapping classes.}
Given a property $p:(\R^d\times\R^D)^m \to \R$ and sets $Q\subseteq \bar Q\subseteq\R^d$, we define families of mappings satisfying $p$ on $Q$ as
\begin{align*}
    \F_p(Q):=\bigg\{F: Q\rightrightarrows \R^D :  &\ p((x_1,F_{x_1}),\ldots ,(x_m,F_{x_m}))\leq 0,\\& \text{for all tuples in $Q^m$, and for all $F_{x_i}\in F(x_i)$, $i\in \llbracket m\rrbracket$}\bigg\}.
\end{align*}
As such, $\F_p(\bar Q)$ and $\F_p(Q)$ describe different objects and cannot be compared (as the mappings in $\F_p(\bar Q)$ and $\F_p(Q)$ are not defined on the same sets). We thus further introduce the family of mappings satisfying $p$ on $\bar Q$, but restricted to $Q$, or equivalently, the mappings $F$ satisfying $p$ on $Q$, and admitting an \emph{extension} to $\bar Q$, i.e. $\exists \bar F\in \F_p(\bar Q)$ such that $\bar F=F$ on $Q$.
\begin{equation*}
    \FQrestr:=\left\{F: Q\rightrightarrows \R^D : \exists \bar F\in \F_p(\bar Q): F(x)=\bar F(x), \ \forall x\in Q\right\}.
\end{equation*}
\paragraph{Interpolation/extension properties.}
Clearly, $\FQrestr \subseteq  \F_p(Q)$ always holds, but the contrary is not necessarily true, as some mappings $F\in\F_p(Q)$ may not admit an extension $\bar F\in \F_p(\bar Q)$. Such mappings thus satisfy $p$ on~$Q$ but cannot be extended to $\bar Q$ in a way as to still satisfy $p$, see, e.g., the example of smooth convex functions \cite[Figure 1]{taylor2017smooth}.

Given a property $p$, we search for a property $\tilde{p}_{\bar Q}$ such that (i) $\F_{\tilde{p}_{\bar Q}}(\bar Q)=\F_p(\bar Q)$, and (ii) mappings satisfying $\tilde{p}_{\bar Q}$ on any $Q\subseteq \bar Q$ always admit an extension in $\F_p(\bar Q)$, i.e., 
\[
    F_{\tilde{p}_{\bar Q}}(Q)=\F_{\tilde{p}_{\bar Q}}(\bar Q)|_Q=\F_{p}(\bar Q)|_Q, \ \forall Q\subseteq \bar Q: \ Q \text{ countable}.
\]
We denote such a condition $\tilde{p}_{\bar Q}$ an \emph{interpolation condition} for $\F_p(\bar Q)$.

\begin{definition}[interpolation/extension condition] \label{def:intrp}
    Given $\bar Q\subseteq \R^d$, and a property $p$, we say $p$ is an interpolation condition for $\F_p(\bar Q)$ if   \begin{equation*} \forall Q\subseteq \bar Q, \ Q \text{ countable}: \ \forall F \in \F_p(Q),\quad \exists F'\in \F_p(\bar Q)\text{ such that } F(x)=F'(x) \  \forall x\in Q. \end{equation*}
\end{definition}
A detailed discussion on differences between this definition and classical ones \cite[Definition 2]{taylor2017smooth} is provided in \cite[Section 2]{rubbens2024constraint}.

\paragraph{Goal and proposed methodology.} Starting from a condition $p$, a domain $\bar Q\subseteq \R^d$ and $K \in \N, K \leq |\bar Q|$, we seek to obtain an interpolation condition for $\F_p(\bar Q)$, which can be $p$ itself. To this end, we construct a sequence of properties $p^{(1)},\, p^{(2)},\, \ldots, p^{(T)}$ such that $\F_{p^{(i)}}(\bar Q)=\F_p(\bar Q)$, but iteratively \emph{stronger} on~$\bar Q$, in the sense that on all $Q\subseteq \bar Q$, $p^{(i+1)}$ is satisfied by less mappings than $p^{(i)}$:  $\F_{p^{(i+1)}}(Q)\subseteq\F_{p^{(i)}}(Q), \ \forall Q \subseteq \bar Q$. Ideally, we hope that the iterative refinement procedure admits a limiting property $p^{(\infty)}$ that will be an interpolation condition.
\subsection{One-point strengthenings}
Consider a property $p:(\R^d\times\R^D)^m \to \R$, a domain $\bar Q \subseteq \R^d$ and some $K\in \mathbb{N},\ m\leq K\leq |\bar Q|$. We take a first step into obtaining an interpolation condition for $\F_p(\bar Q)$ and build from $p$ a stronger definition of this class, which ensures that any mapping $F\in \F_p(Q')$ where $Q'\subseteq \bar Q$ and $|Q'|=K$, admits an extension satisfying $p$ to any $z\in \bar Q$. Requiring an extension to any $z\in \bar Q$ individually is weaker than requiring an extension to $\bar Q$, i.e., an extension to any $z\in \bar Q$ simultaneously. 

We denote by $\tilde{p}_{K,\bar Q}:(\R^d\times\R^D)^K \to \R$ such a stronger definition, that thus satisfies (i) $\F_{\tilde{p}_{K,\bar Q}}(\bar Q)=\F_p(\bar Q)$, (ii) $\forall Q\subseteq \bar Q$, $\F_{\tilde{p}_{K,\bar Q}}(Q)\subseteq \F_p(Q)$, and (iii) $\forall Q'\subseteq\bar Q$: $|Q'|=K$, and for any $z\in \bar Q$, $F\in \F_{\tilde{p}_{K,\bar Q}}(Q')\Leftrightarrow F\in \F_p(Q'\cup z)|_{Q'}$. Observe that the number of points and mapping values taken as arguments by $p$ and $\tilde{p}_{K,\bar Q}$, respectively, varies from $m$ to $K$.

We construct $\tilde{p}_{K,\bar Q}$ as follows. Given $Q'=\{(x_1,\ldots,x_K)\} \subseteq \bar Q, \ |Q'|=K$, $F\in \F_p(Q')$, and $z\in \bar Q$, it is possible to extend $F$ from $Q'$ to $z$ if there exists some $F_z\in \R^D$ compatible with $p$ and $F$ on $Q'$, i.e., if and only if
\begin{align*}
    0\geq \min_{\substack{\tau \in \R,\\F_z\in \R^D}}\ \ \{\tau \text{ s.t. } &p(\pi((x_1,F_{x_1}),\ldots ,(x_{m-1},F_{x_{m-1}}),(z,F_z)))\leq \tau,\\& \forall \text{ tuple in $Q^{'m-1}$, $F_{x_i}\in F(x_i)$, $i\in\llbracket m-1\rrbracket$, and permutation $\pi$}\}.
\end{align*}
Since this must hold for any $z\in \bar Q$, $\tilde{p}_{K,\bar Q}$ should enforce it for the ``worst'' $z$, i.e.,
\begin{align}\label{eq:def1}
    0\geq \underset{z\in \bar Q}{\max} \  \min_{\substack{\tau \in \R,\\F_z\in \R^D}}\ \{\tau \text{ s.t. } &p(\pi((x_1,F_{x_1}),\ldots ,(x_{m-1},F_{x_{m-1}}),(z,F_z)))\leq \tau,\\& \forall \text{ tuple in $Q^{'m-1}$, $F_{x_i}\in F(x_i)$, $i\in\llbracket m-1\rrbracket$, and permutation $\pi$}\}.\nonumber
\end{align}
Letting $\tilde{p}_{K,\bar Q}$ be the solution to the optimization problem defined in \eqref{eq:def1} implies that any $F\in \F_{\tilde{p}_{K,\bar Q}}(Q')$, hence satisfying $\tilde{p}_{K,\bar Q}((x_1,F_{x_1}),\ldots ,(x_{K},F_{x_{K}}))\leq 0$, admits an extension in $\F_p(Q'\cup z)$, for any $z\in \bar Q$. In particular, letting $z\in Q'$ naturally implies $F\in \F_p(Q')$.

We now formally define $\tilde{p}_{K,\bar Q}$ and prove that $\F_{\tilde{p}_{K,\bar Q}}(\bar Q)=\F_p(\bar Q)$, and $\tilde{p}_{K,\bar Q}$ is stronger than $p$ on~$\bar Q$, before illustrating the definition on the example of smooth convex functions.
\begin{definition}[$K$-wise one-point strengthening] \label{def:1pt_strengthening}
    Consider a property $p:(\R^d\times\R^D)^m \to \R$, $ \bar Q \subseteq \R^d$,  and $K\in \mathbb{N}, m\leq K\leq |\bar Q|$. We define $\tilde{p}_{K,\bar Q}:(\R^d\times\R^D)^K \to \R$, the $K$-wise one-point strengthening of $p$, as
    \begin{equation}\label{eq:1pt_strengthening}
        \begin{aligned}[t]
            \tilde{p}_{K,\bar Q}((x_1,F_{x_1}),\ldots ,(x_{K},F_{x_{K}}))=\underset{z\in \bar Q}{\max} \  \underset{F_z\in \R^D,\tau \in \R}{\min}\ \{\tau \text{ s.t. } &p(\pi((x_1,F_{x_1}),\ldots ,(x_{m-1},F_{x_{m-1}}),(z,F_z)))\leq \tau,\nonumber\\& \forall \text{ tuple in $\{x_1,\ldots,x_K\}^{m-1}$, permutation $\pi$},\\& \text{and $\forall F_{x_i}\in F(x_i)$, $i\in\llbracket m-1\rrbracket$}\}.
        \end{aligned}
    \end{equation}
\end{definition}
\begin{lemma}
    Given $\bar Q \subseteq \R^d$, $K \in \N, K\leq |\bar Q|$, consider a property $p$ and its one-point strengthening $\tilde{p}_{K,\bar Q}$, as defined in \eqref{eq:1pt_strengthening}.
    Then, $\forall Q\subseteq \bar Q$, $\F_{\tilde{p}_{K,\bar Q}} (Q)\subseteq \F_p(Q)$, and $\F_p(\bar Q)=\F_{\tilde{p}_{K,\bar Q}} (\bar Q)$.
\end{lemma}
\begin{proof}
    Taking $z\in Q$ in \eqref{eq:1pt_strengthening}, it holds that all $F \in \F_{\tilde{p}_{K,\bar Q}} (Q)$ satisfy $p$ on~$Q$, hence the first implication holds. Similarly, all $F \in \F_{\tilde{p}_{K,\bar Q}}(\bar Q)$ satisfy $p$ on $\bar Q$, hence $\F_{\tilde{p}_{K,\bar Q}}(\bar Q)\subseteq \F_{p}(\bar Q)$.

    Suppose now by contradiction the existence of some $F \in \F_p(\bar Q)$ which does not belong to $\F_{\tilde{p}_{K,\bar Q}}(\bar Q)$. Then, at some $(x_1,\ldots,x_{m-1}) \in \bar Q$, for some $F_{x_i}\in F(x_i)$, $i\in\llbracket m-1\rrbracket$, and for some permutation $\pi$, it holds that
    \begin{equation*}
        0< \max_{z\in \bar Q}\min_{\substack{\tau \in \R, \ F_z \in \R^D}} \  \tau \text{ s.t. }  p(\pi((x_1,F_{x_1}),...,(x_{m-1},F_{x_{m-1}},(z,F_z)))\leq \tau.
    \end{equation*}
    Hence, at some $z\in \bar Q$ there exists no $F_z$ extending $p$ from $(x_1,\ldots,x_{m-1})$ to $z$, in contradiction with the assumption $F\in \F_p(\bar Q)$.
\end{proof}
\begin{remark}\label{rem:param}
    \Cref{def:1pt_strengthening} allows evaluating $\tilde{p}_{K,\bar Q}$ on any set $Q\subseteq \bar Q$, by simply enforcing the solution to \eqref{eq:1pt_strengthening} to be non-positive for all tuples in $Q^K$. Alternatively, one can obtain a closed-form expression of $\tilde{p}_{K,\bar Q}$ by computing the solution to \eqref{eq:1pt_strengthening} in a \emph{parametric} way in $\{x_1,\ldots,x_K\}$, i.e., keeping these as parameters. Evaluating $\tilde{p}_{K,\bar Q}$ on any set $Q\subseteq \bar Q$ then reduces to simply imposing a set of conditions on $F:Q\to \R^D$, without having to solve an optimization problem.
\end{remark}

In what follows, when $K,\bar Q$ are clear from the context, we drop the explicit dependence on $K,\bar Q$ and denote $\tilde p_{k,\bar Q}$ via the simpler $\tilde p$ instead.
\begin{example}[smooth convex functions]
    Let $Q=\{(x_1,x_2)\}$ and $F:Q\to\R^{1+d}:x_i\to (f_i,g_i)$, where $g_i$ is a priori unrelated from $f_i$ while meant to belong to its differential, see \cref{rem:diff}. The $2$-wise (pairwise) one-point strengthening of \eqref{eq:convexity_quadratic_UB} with respect to $\R^d$ is given by
    \begin{align}\label{eq:quadratic_UB_1pt}\tag{$\tilde p'_{0,L}$}
        \tilde p'_{0,L}((x_1,f_1,g_1),(x_2,f_2,g_2)) =         \max_{z\in \R^d} \min_{\substack{\tau \in \R, \\ g_z\in \R^d,\ f_z\in \R}} \  \tau \text{ s.t.  }  &f_z- f_i-\langle g_i,z-x_i\rangle -\frac{L}{2}\|z-x_i\|^2\leq \tau,\ && i=1,2 \nonumber\\
        &f_i-f_z+\langle g_z,z-x_i\rangle -\frac{L}{2}\|z-x_i\|^2\leq \tau,\ &&i=1,2\nonumber\\
        &f_z- f_i-\langle g_z,z-x_i\rangle \leq \tau,\ &&i=1,2\nonumber\\
        &f_i- f_z+\langle g_i,z-x_i\rangle \leq \tau,\ &&i=1,2.\nonumber
    \end{align}
\end{example}

\paragraph{Stability of a strengthening.}
Given $\bar Q\subseteq\R^d$, $K\in \N,K\leq |\bar Q|$, and starting from a given property $p$, one can strengthen $p$ by iteratively computing one-point strengthenings $p^{(1)}=\tilde{p}_{K,\bar Q},\  p^{(2)}=\tilde{p}_{K,\bar Q}^{(1)},\ldots$, with the goal of reaching a property $p_{K,\bar Q}^{(\infty)}$ that is \emph{stable under $K$-wise one-point strengthening}, i.e., $p_{K,\bar Q}^{(\infty)} =\tilde{p}_{K,\bar Q}^{(\infty)}$ in the sense that $\F_{p_{K,\bar Q}^{(\infty)}}(Q)=\F_{\tilde{p}_{K,\bar Q}^{(\infty)}}(Q)$, $\forall Q\subseteq Q'$.

\begin{definition}[stability under one-point strengthening]
    Given $\bar Q\subseteq\R^d$, $K\in \N,K\leq |\bar Q|$, and a property $p$, we say $p$ is stable under $K$-wise one-point strengthening with respect to $\bar Q$ if $\tilde{p}_{K,\bar Q}=\ p$, where $\tilde{p}_{K,\bar Q}$ is defined in \eqref{eq:1pt_strengthening}.
\end{definition}
The iterative procedure can, for instance, reach such a stable property in finite time, see, e.g., \cref{prop:1pt_strengthening_smooth_convex_UB}, or as a limit to which the sequence of properties, under their closed-form, converges, see, e.g., \cref{prop:1pt_strengthening_smooth_convex}.
\paragraph{$K$-wise one-point extensible properties and extension properties.}
The final property of the refinement procedure, $p^{(\infty)}$, is $K$-wise one-point extensible, in the sense that it can be extended from any subset $Q\subseteq \bar Q$ of cardinality $K$, to any $z \in \bar Q$. To obtain a \emph{global one-point extensible property} for $\bar Q$, i.e., a property that can be extended from any finite subset to any additional point, it then suffices to relaunch the procedure for a sequence of larger and larger $K$, i.e., $K_1\leq K_2\leq \ldots\leq |\bar Q|$, with typically $K_i=K_{i-1}+1$, to obtain a sequence $p^{(\infty)}_{K_1,\bar Q},p^{(\infty)}_{K_2,\bar Q},\ldots,p^{(\infty)}_{K_k,\bar Q}$, where each $p^{(\infty)}_{K_k,\bar Q}$ is computed by initializing the procedure at $p^{(0)}_{K_i,\bar Q}=p^{(\infty)}_{K_{i-1},\bar Q}$.

In addition, under reasonable regularity assumptions, interpolation conditions and global one-point extensible conditions are equivalent~\cite[Theorem 2] {rubbens2024constraint}.
\subsection{Computing of one-point strengthenings }\label{sec:practical_comments}
Starting from a given $p$, obtaining $\tilde{p}_{K,\bar Q}$ requires the ability to solve~\eqref{eq:1pt_strengthening}. We tackle this task by dualizing the inner problem in \eqref{eq:1pt_strengthening}, which can be done explicitly for a range of natural conditions used in the optimization community. Afterwards, whenever there is no duality gap, solving \eqref{eq:1pt_strengthening} amounts to solving a single maximization problem, whose number of variables grows with the number of inequalities in $p$, and $K$. There is no general recipe to solve this maximization problem, which often admits several solutions that can lose the initial structure of $p$ in $F$, e.g., linearity.

Hence, in most cases, we do not solve~\eqref{eq:1pt_strengthening} exactly but instead propose relaxations of $\tilde{p}_{K,\bar Q}$, which might not be optimal solutions to \eqref{eq:1pt_strengthening}, but are nevertheless stronger properties than $p$.
In particular, we often conserve only a limited number of non-zero dual coefficients for the inner problem. Then, when the resulting conditions are polynomials in the remaining coefficients, we reformulate/relax them as SDPs, see \cref{sec:related_works}. Alternatively, one can choose to impose particular numerical values for the non-zero coefficients. Both relaxations of $\tilde{p}_{K,\bar Q}$ conserve the same structure in terms of $F$ as $p$ and can thus be straightforwardly relied upon when imposed on datasets that are variables to given optimization problems.
\begin{remark}\label{rem:complexity}
    Note that relaxing $\tilde{p}_{K,\bar Q}$ involves trading off between (i)~obtaining conditions as close as possible to $\tilde{p}_{K,\bar Q}$, e.g., via convex reformulations (e.g., co-positive reformulations, semi-definite, etc.), and (ii)~obtaining conditions whose one-point strengthening are easily computable.
\end{remark}

For similar reasons, we are, in general, not able to prove stability of a one-point strengthening. For exploring this venue, we provide some insight into the iterative process on two examples for which extension conditions are known: that of smooth convex functions (see \cref{sec:smooth_convex}) and that of Lipschitz operators (see \cref{app:op_Lip}).

\subsection{Operator and function classes}
\Cref{sec:formal_extension} introduced the proposed iterative refinement procedure in a general context. For the sake of readability, we now introduce lighter notation related to the cases we analyze in what follows.

\paragraph{Function classes and differentially consistent properties.}\label{rem:diff}
We first consider function classes characterized by properties involving function values and first-order derivatives, e.g., the class of $L$-smooth convex functions and its definition~\eqref{eq:true_smooth_convexity}. Mappings of interest then take the form $F:Q \rightrightarrows \R^{1+d}:\ x\rightrightarrows (f(x),g(x))\in \R \times \R^d$. While meant to satisfy $g(x)\in \partial f(x)$, the mappings $g(x)$ and $f(x)$ are \emph{a priori} unrelated, other than satisfying a certain algebraic relationship. We call a property $p$ whose satisfaction everywhere implies
$g(x)\in \partial f(x),\ \forall x\in \R^d$ \emph{differentially consistent} \cite[Definition 2]{rubbens2024constraint}.  We denote by $S=\{(x_i,f_i,g_i)\}_{i \in \mathcal{I}}$ for an index set $\mathcal{I}$ the set of points $Q$ together with the associated value of mapping $F\in \F_p(Q)$ to be extended, where a multivalued mapping is accounted for by letting several points in $S$ be equal, with different associated values of $g_i$'s.
\paragraph{Compact notation for a property and its strengthening.} For any property $p$, we denote (i)~with an index (e.g., $p_\alpha$) the possible parameters for this property, (ii)~with a superscript the points at which it is evaluated (e.g., $p_\alpha^{xy}=p_\alpha((x,f(x),g(x)),(y,f(y),g(y))$ or $p_\alpha^{ij}=p_\alpha((x_i,f_i,g_i),(x_j,f_j,g_j))$), (iii)~with a tilde the corresponding one-point strengthening~\eqref{eq:1pt_strengthening} of $p$ (e.g., $\tilde{p}_\alpha$ the one-point strengthening of $p_\alpha$ obtained by~\eqref{eq:1pt_strengthening}), and (iv)~with a superscript in between parentheses the iterates of the one-point strengthening procedure (e.g., $p_\alpha^{(i+1)}=\tilde{p}_\alpha^{(i)}$) -- we never use two types of superscripts simultaneously.

\paragraph{Operator classes.} \cref{sec:operators} and \cref{sec:operators2} provide examples of property strengthenings for the case of operators (as opposed to functions). Mappings of interest take the form $F:\R^d \to \R^d: x\to t(x)$, and in this case we denote by $S=\{(x_i,t_i)\}_{i \in \mathcal{I}}$ the dataset to be extended, where again a multivalued mapping is accounted for by letting several points in $S$ be equal, with different associated $t_i$'s. To avoid notation clashes, we refer to properties of operators using the letter $q$ instead of $p$, with all variations around this notation also used (e.g., $q_\alpha$ and its one-point strengthening $\tilde{q}_\alpha$).

Finally, we denote $2$-wise by pairwise, $3$-wise by terwise, etc.
\subsection{Example: smooth convex functions} \label{sec:smooth_convex}
\begin{figure}
    \centering
    \includegraphics[width=0.5\linewidth]{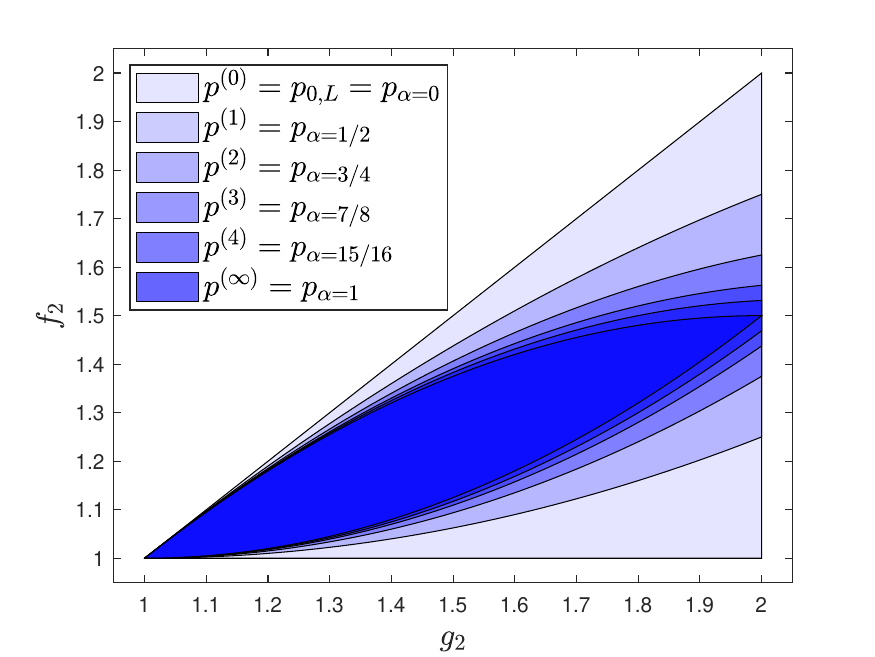}
    \caption{From \cite[Figure 1]{rubbens2024constraint}. Given $L=1$, $x_1=f_1=0, \ g_1=x_2=1$, allowed regions for $f_2$ as a function of $g_2$, according to $p^{(i)}=p_{\alpha_i}$, $i=0,\ldots,4$, where $\alpha_0=0$, $\alpha_{i+1}=(1+\alpha_i)/2$, and $p_{\alpha_i}$ is defined in \eqref{eq:initial_condition}, and according to $p^{(\infty)}=p_{\alpha=1}$=\eqref{eq:true_smooth_convexity}. }
    \label{fig:smooth_convex}
\end{figure}

This section illustrates the refinement procedure on a known and well-studied case: that of smooth convex functions (\cite{nesterov2018lectures}). We argue that it is a perfect example to illustrate how the technique can be used, as interpolation conditions \eqref{eq:true_smooth_convexity} are known for this class.
We show how to obtain ~\eqref{eq:true_smooth_convexity}, first from \eqref{eq:convexity_quadratic_UB} and then from the most natural definition \eqref{eq:convexity_gradient_lipschitz}, without resorting to integration -- which is the classical way to obtain those inequalities, albeit arguably not a natural argument given the nature of the problem.
To this end, consider first the following one-point strengthening.
\begin{proposition} \label{prop:1pt_strengthening_smooth_convex_UB}
    Let $Q=\{(x_1,x_2)\}$ and consider \eqref{eq:convexity_quadratic_UB}. Denote by $\tilde p'_{0,L}$ the pairwise one-point strengthening of \eqref{eq:convexity_quadratic_UB} on $\R^d$. It holds that $\tilde p'_{0,L}$=\eqref{eq:true_smooth_convexity}.
\end{proposition}
\begin{proof}
    The one-point strengthening of \eqref{eq:convexity_quadratic_UB}, \eqref{eq:quadratic_UB_1pt}, can be relaxed to
    \begin{align*}
        \tilde p_{0,L}^{'ij} \geq         \max_{z\in \R^d}\min_{\substack{\tau \in \R,\\ g_z \in \R^d, \ f_z \in \R}} \ \tau \text{ s.t.  }  &f_z- f_i-\langle g_i,z-x_i\rangle -\frac{L}{2}\|z-x_i\|^2\leq \tau, \ i=1,2&& (\lambda_i)\nonumber\\
        &f_j- f_z+\langle g_j,z-x_j\rangle \leq \tau, \ j=1,2,&& (\mu_j)\nonumber
    \end{align*}
    where $\lambda_i$ and $\mu_j$ are the associated dual coefficients. Setting $\lambda_j=\mu_i=0, \ \lambda_i=\mu_j=\frac{1}{2}$ for some choice $i\neq j=1,2$, and dualizing the inner problem yields
    \begin{align*}
        \tilde p_{0,L}^{'ij} \geq \frac{1}{2}\max_{z\in \R^d} \ f_j+\langle g_j,z-x_j\rangle- f_i-\langle g_i,z-x_i\rangle -\frac{L}{2}\|z-x_i\|^2, \textcolor{white}{blblablablabablabllal}
    \end{align*}
    which is maximized in $z=x_i+\frac{g_j-g_i}{L}$, hence
    \begin{align*}
        0\geq \tilde p_{0,L}^{'ij} \Rightarrow 0 \geq f_j-f_i+\langle g_j,x_i-x_j\rangle+\frac{1}{2L}\|g_j-g_i\|^2\Leftrightarrow
        0 \geq p_{0,L}^{''ij}.\textcolor{white}{blalbblabll}
    \end{align*}
    Equivalency between $\tilde p'_{0,L}$ and \eqref{eq:true_smooth_convexity} follows from the knowledge that the latter is an interpolation condition, hence cannot be further reinforced.
\end{proof}
Starting from \eqref{eq:convexity_quadratic_UB}, we thus recover \eqref{eq:true_smooth_convexity} after a single round of one-point strengthenings. As a second step in our example, we show how to recover \eqref{eq:true_smooth_convexity} from \eqref{eq:convexity_gradient_lipschitz}, which requires several iterations of one-point strengthenings that all share the same structure. \Cref{fig:smooth_convex} illustrates the iterative refinement procedure for $\F_{0,L}$.  We first present a single iteration of the strengthening procedure starting with~\eqref{eq:convexity_gradient_lipschitz} (when $\alpha=0$).

\begin{proposition} \label{prop:1pt_strengthening_smooth_convex}
    Consider $Q=\{(x_1,x_2)\}$, $\alpha\geq 0$, and $p_\alpha$ a pairwise condition defined by $\forall i,j =1,2$
    \begin{align} \tag{$p_\alpha$}
        p^{ij}_\alpha\leq 0 \Leftrightarrow\begin{cases}
            &f_i\geq f_j+\langle g_j,x_i-x_j\rangle+\frac{\alpha}{2L}\|g_i-g_j\|^2\ \label{eq:initial_condition}
            \\
            &\|g_i-g_j\|^2\leq L^2 \|x_i-x_j\|^2.
        \end{cases}\textcolor{white}{blabl}
    \end{align}
    Then, $\tilde{p}_{\alpha}$, the pairwise one-point strengthening of \eqref{eq:initial_condition} with respect to $\R^d$, is satisfied only if, $\forall i,j =1,2$
    \begin{align}
        \tilde p_\alpha^{ij}\leq 0 \Rightarrow\begin{cases}
            &f_i\geq f_j+\langle g_j,x_i-x_j\rangle+\frac{(1+\alpha)/2}{2L}\|g_i-g_j\|^2\ \label{eq:final_condition}
            \\
            &\|g_i-g_j\|^2\leq L^2 \|x_i-x_j\|^2.
        \end{cases}
    \end{align}
\end{proposition}
\begin{proof}
    The one-point strengthening of \eqref{eq:initial_condition} is given by
    \begin{align*}
        \tilde{p}_{\alpha}(f,Q)= \max_{z\in \R^d}\min_{\substack{\tau \in \R, \\ g_z\in \R^d,\ f_z\in \R}} \tau
        \text{ s.t. }&
        -f_z+ f_i+\langle g_i,z-x_i\rangle+ \frac{\alpha}{2L}\|g_i-g_z\|^2\leq \tau, \ i=1,2 && (\lambda_i)\nonumber\\
        & -f_i+ f_z-\langle g_z,z-x_i\rangle+  \frac{\alpha}{2L}\|g_i-g_z\|^2\leq\tau , \ i=1,2&&(\mu_i)\nonumber\\
        & \frac{1}{2L}\|g_z-g_i\|^2-\frac{L}{2}\|z-x_i\|^2\leq \tau, \ i=1,2, &&(\delta_i)\nonumber
    \end{align*}
    where $\lambda_i, \mu_i, \delta_i, \ i=1,2,$ are the associated dual coefficients. Partial dualization of the inner problem over $\tau$ and $f_z$, and minimization over $g_z$ yields
    \begin{align*}
        \tilde{p}_{\alpha}(f,Q)\geq \max_{\substack{z\in \R^d,\\ \lambda_i, \ \mu_i, \ \delta_i\geq 0\\\sum_{i} (\lambda_i+\mu_i+ \delta_i)=1}} \ &\sum_{i=1}^2 \big(\lambda_i (f_i+\langle g_i,z-x_i\rangle+ \frac{\alpha}{2L}\|g_i-g_z^\star\|^2)  \\&\quad+ \mu_i (-f_i-\langle g_z^\star,z-x_i\rangle+  \frac{\alpha}{2L}\|g_i-g_z^\star\|^2 )+\frac{\delta_i}{2L} (\|g_z^\star-g_i\|^2-L^2\|z-x_i\|^2)\big)\nonumber \\
        \text{ s.t. }
        0&=\sum_{i=1}^2 (\lambda_i-\mu_i), \nonumber\\
        0& =\sum_{i=1}^2\left((\mu_i+\lambda_i) \frac{\alpha}{L} (g_z^\star-g_i)-\mu_i (z-x_i)+\frac{\delta_i}{L} (g_z^\star-g_i)\right)\nonumber,
    \end{align*}
    where $g_z^\star$ minimizes the inner problem in $\tilde{p}_{\alpha}$.

    In particular, setting $\lambda_i=0, \mu_i=0, \delta_i=\frac{1}{2}, \  i=1,2$, $g_z^\star=\frac{g_1+g_2}{2}$ and $z=\frac{x_1+x_2}{2}$ yields a feasible solution of $\|g_1-g_2\|^2-L^2\|x_1-x_2\|^2$, hence $\tilde{p}_{\alpha}$ cannot be satisfied unless this quantity is negative. In addition, setting  $\lambda_i=0, \  \mu_j=0, \ \delta_j=0, \ \lambda_j=\mu_i=\delta_i=\frac{1}{3}$, $z=x_i+\frac{g_j-g_i}{2L}$ and $g_z^\star=\frac{g_i+g_j}{2}$ yields a feasible solution of $-f_j+ f_i+\langle g_i,x_j-x_i\rangle+ \frac{1+\alpha}{2} \frac{1}{2L}\|g_i-g_j\|^2$, $i,j=1,2$. Hence, $\tilde{p}_{\alpha}$ is at least as strong as \eqref{eq:final_condition}.
\end{proof}
To recover \eqref{eq:true_smooth_convexity} from \eqref{eq:convexity_gradient_lipschitz}, it suffices to successively apply the one-point strengthening of \cref{prop:1pt_strengthening_smooth_convex}, starting with $p^{(0)}= $\eqref{eq:convexity_gradient_lipschitz}, to reach $p^{(\infty)}=$ \eqref{eq:true_smooth_convexity}.
\begin{proposition}
    Consider $Q=\{(x_1,x_2)\}$, \eqref{eq:convexity_gradient_lipschitz}, the sequence $\{\alpha_k\}_{k=0,\ldots}$, initiated at $\alpha_0=0$ and defined by
    \begin{equation}
        \alpha_{k+1}=\frac{\alpha_k+1}{2},\label{eq:sequence}
    \end{equation}
    and the associated sequence of conditions $\{p^{(k)}\}_{k=0,\ldots}:=\{p_{\alpha_k}\}_{k=0,\ldots}$, where $p_{\alpha_k}$ is defined in \eqref{eq:initial_condition}.

    Then, $p^{(k+1)}$ is a relaxation of $\tilde{p}^{(k)}$ as defined in \eqref{eq:final_condition}. In addition, $p^{(\infty)}$,  the final result of the refinement procedure defined in \cref{sec:formal_extension}, satisfies $p^{(\infty)}=\tilde{p}_{K,\bar Q}^{(\infty)}=p_{\alpha_{(\infty)}}=p_{\alpha=1}$=\eqref{eq:true_smooth_convexity}.
\end{proposition}
\begin{proof}
    By \cref{prop:1pt_strengthening_smooth_convex} and definition of $\{(\alpha_k)\}_{k=0,\ldots,}$, the first statement holds.

    In addition, $\{(\alpha_k)\}_{k=0,\ldots,}$ converges to $\alpha_{(\infty)}=1$. Hence, applying recursively \cref{prop:1pt_strengthening_smooth_convex} starting from $\alpha=0$ allows iteratively strengthening $p$ until $p^{(\infty)}=p_{\alpha_{(\infty)}}=$\eqref{eq:true_smooth_convexity}, which is stable under one-point strengthening.
\end{proof}

\paragraph{Concluding observations.} In this example, the iterative strengthening procedure described in \cref{sec:formal_extension} converges towards the same property even when initialized from two different descriptions of $\F_{0,L}(\R^d)$. This holds despite relaxations in the computations of the one-point strengthenings involved. In addition, in both cases, the iterative procedure even converges towards an interpolation condition for $\F_{0,L}(\R^d)$. However, the initial condition given to the procedure has a non-trivial impact on the number of iterations required to reach this condition, which varies from a single iteration to an asymptotic convergence. Note that the same observations apply to the other known case on which we applied the procedure: that of Lipschitz operators, presented in \cref{app:op_Lip}.

\section{Applications: improved descriptions of a few problem classes}\label{sec:applications}
This section strengthens the description of a few classical classes of functions or operators. For all of them, we perform a single round of the refinement procedure described in \cref{sec:formal_extension}. We propose several reformulations/relaxations of these one-point strengthenings and show that the corresponding refined descriptions are strictly stronger than the classical descriptions serving as starting points.

A particularly common scenario is that of initial properties that are linear in the function values and the Gram matrix of the point and subgradient (or operator) values. That is, properties that are linear in $f_i$'s as well as linear in terms of $\langle g_i,x_i\rangle,\ \|g_i\|^2,\ \|x_i\|^2$ (or $t_i$ instead of $g_i$ for operators). In this case, we will see that the resulting strengthening can be expressed as requiring non-negativity of a degree $2k$ polynomial (for some $k\in\mathbb{N}$), whose coefficients are linear in $f_i$'s and the Gram matrix. While the set of non-negative polynomials is generally intractable, it can be approximated, e.g., by constraining the polynomial to lie within the set of SOS polynomials \cite{hilbert1888darstellung,lasserre2009moments,marshall2008positive}, which can be expressed as an SDP, linear in $f_i$'s and Gram matrix, see \cref{sec:related_works} for details. In particular, for one-dimensional polynomials, or two-dimensional polynomials when $k=2$, SOS-polynomials and non-negative polynomials are strictly equivalent. \cref{tab:summary} summarizes these applications of the refinement procedure.
\begin{table}[t!]
    \centering
    \begin{tabular}{@{}lll@{}}
        \toprule
        \textbf{Class}& \textbf{Section} & \textbf{Proposed formulations/relaxations}  \\
        &  & \textbf{of the one-point strengthening}  \\
        \midrule
        Smooth functions satisfying & \ref{sec:loja} & \eqref{loja_strengthened_1}: Pairwise linear inequalities, involving
        \\
        a \L{}ojasiewicz condition & & \textcolor{white}{(3.1):} maximization over a single variable in $\R^d$,
        \\
        &  &\eqref{def:matrix_loja}: SDP exact reformulation of \eqref{loja_strengthened_1}.\\
        \midrule
        Convex blockwise smooth & \ref{sec:blockwise}& \eqref{eq:blockwise_strengthened}: Pairwise linear inequalities, involving
        \\
        functions& & \textcolor{white}{(3.3):} maximization over a single variable in $\R^d$, \\

        &&\eqref{def:matrix_block}: SDP exact reformulation of \eqref{eq:blockwise_strengthened},\\
        &  & \eqref{eq:block_quartic}: Terwise quadratic inequalities, involving a $if$ condition.\\
        \midrule
        Lipschitz strongly & \ref{sec:operators} & \eqref{eq:final_p_op}: Terwise linear inequalities, involving \\
        monotone operators & & \textcolor{white}{(3.3):} maximization over 3 variables in $\R^d$,\\&&\eqref{def:matrix_op}: SDP relaxation of \eqref{eq:final_p_op}.\\ \midrule
        Strongly monotone & \ref{sec:operators2} &
        \eqref{eq:final_popi}:  Terwise linear inequalities, involving\\
        cocoercive operators & & maximization over 3 variables in $\R^d$,\\&& \eqref{def:matrix_opi}: SDP relaxation of \eqref{eq:final_popi}.\\
        \midrule
        Uniformly convex functions & \ref{sec:unif_conv} & \eqref{unif_conv_strengthened}: Terwise convex inequalities, involving \\
        && maximization over a single variable in $\R^d$ .\\ \midrule
        Smooth convex functions& \ref{sec:extension_subset}
        & Pairwise linear inequalities. \\
        on convex subsets of $\R^d$&& \\
        \bottomrule
    \end{tabular}
    \caption{Summary: function and operator classes studied within this work, along with the nature of the strengthened conditions as expressions of the function values and Gram matrix (of $(x_i,g_i)$ or~$(x_i,t_i)$).}
    \label{tab:summary}
\end{table}
\subsection{Extension of smooth functions satisfying a \L{}ojasiewicz condition}\label{sec:loja}
This section studies the class $\F_{\text{\L{}}_{\mu,L}}(\R^d)$ of $L$-smooth functions $f$ satisfying a quadratic \L{}ojasiewicz condition, that is satisfying, $\forall x,y \in \R^d$ \cite{abbaszadehpeivasti2022conditionslinearconvergencegradient,bolte2007lojasiewicz,bolte2010characterizations,lojasiewicz1963propriete},
\begin{equation}\tag{$\pL$}
    \pL^{xy}\leq 0\Leftrightarrow \left\{\begin{aligned}
            f(x)&\leq f_\star+\frac{\|\nabla f(x)\|^2}{2 \mu } \label{eq:loja_smooth}
            \\
            f(x)&\geq f_\star\\
            f(x)&\geq f(y)+\frac{1}{2}\langle \nabla f(x)+\nabla f(y), x-y\rangle+\frac{1}{4L}\|\nabla f(x)-\nabla f(y)\|^2-\frac{L}{4} \|x-y\|^2,
    \end{aligned}\right.
\end{equation}
where $f_\star=\min_x f(x)$ is guaranteed to exist and $0\leq \mu \leq L$. The last condition imposes smoothness of the function, equivalently to $|f(x)-f(y)-\langle \nabla f(y), x-y\rangle|\leq \frac{L}{2}\|x-y\|^2$ \cite[Theorem 4]{taylor2017exact}.

\begin{figure}
    \vspace{-0.3cm}
    \centering
    \includegraphics[width=0.5\linewidth]{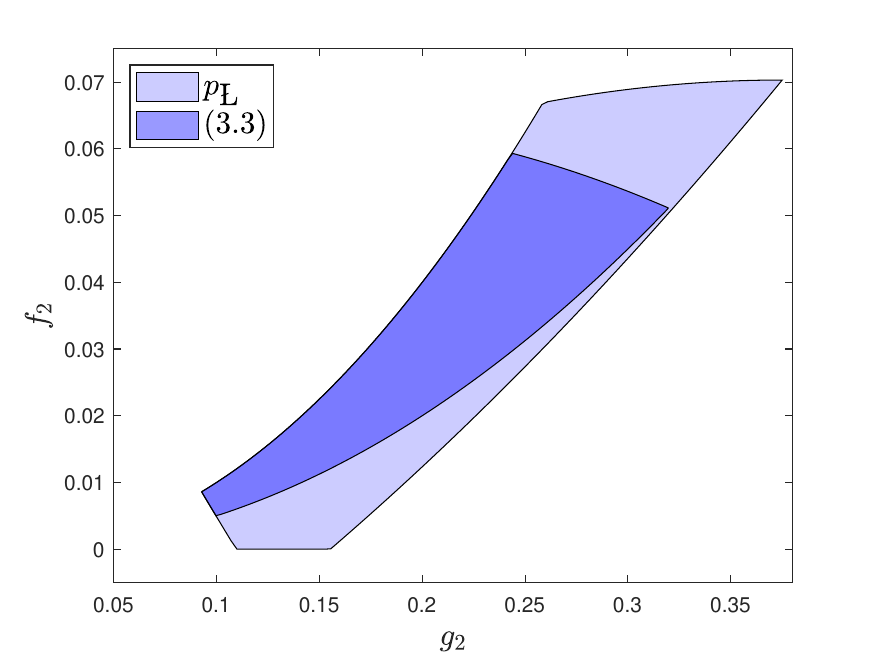}
    \caption{Allowed regions for $f_2$ as a function of $g_2$, according to \eqref{loja_strengthened_1} and \eqref{eq:loja_smooth} given $\mu=1/2, \ L=1$, $(x_\star,f_\star,g_\star)=(0,0,0), \ (x_1,f_1,g_1)=(1, \frac{1}{4}, \frac{1}{2})$ and $x_2=\frac{3}{8}$.}
    \label{fig:loja}
\end{figure}

Exploiting the iterative refinement technique, \eqref{eq:loja_smooth} can be strengthened as follows.
\begin{proposition} \label{prop:1pt_strengthening_loja_smooth}
    Consider $Q=\{(x_1,x_2,x_\star)\}$, and \eqref{eq:loja_smooth}. Then, $\tpL$, the pairwise one-point strengthening of \eqref{eq:loja_smooth} with respect to $\R^d$, is satisfied only if
    \begin{equation}
        \tpL^{ij}\leq 0 \Rightarrow \left\{
            \begin{aligned}
                f_i&\leq f_\star+\frac{\|g_i\|^2}{2 \mu }
                \\
                f_i&\geq f_\star+\frac{\|g_i\|^2}{2 L}\\
                f_i&\geq f_j+\frac{1}{2}\langle g_i+g_j, x_i-x_j\rangle+\frac{1}{4L}\|g_i-g_j\|^2-\frac{L}{4} \|x_i-x_j\|^2
                \\\MoveEqLeft[-1]+\underset{0\leq \alpha \leq  \frac{2\mu}{2L+\mu}}{\max } \frac{\alpha}{(1-\alpha)(2\mu-(L+\mu)\alpha)} \\
                \MoveEqLeft[-5.5]\cdot \left((1-\alpha)^2(L+\mu)(f_i-f_\star-\frac{\|g_i\|^2}{2L})-(L-\mu)(f_j-f_\star+\frac{\|g_j\|^2}{2L})\right).
        \end{aligned}\right.
        \label{loja_strengthened_1}
    \end{equation}
\end{proposition}
The detailed proof of \cref{prop:1pt_strengthening_loja_smooth} is deferred to \cref{app:loja_smooth}. It consists in (i) formulating the pairwise one-point strengthening \eqref{eq:1pt_strengthening} of $\tpL$, (ii) dualizing the inner problem (minimization over $\tau \in \R$, $f_z \in \R$, $g_z\in \R^d$) in this one-point strengthening as to relax it into a maximization problem over positive dual variables and $z\in \R^d$, (iii) explicitly maximizing this relaxation in $z$, and (iv) computing several feasible solutions to this relaxation, by setting some dual variables to $0$ and either explicitly maximize over the remaining variables (e.g., to obtain $f_i\geq f_\star+\frac{\|g_i\|^2}{2L}$), or obtain a constraint involving one of the remaining variables (e.g., $\alpha$ in the last condition in \eqref{loja_strengthened_1}).
\begin{remark}
    Setting $\alpha=0$ in the last condition in \eqref{loja_strengthened_1} allows recovering the last condition in \eqref{eq:loja_smooth}, hence these conditions, although relaxations of $\tpL$, are indeed at least as strong as \eqref{eq:loja_smooth}. The example illustrated in \cref{fig:loja} allows concluding \eqref{loja_strengthened_1} is strictly stronger than \eqref{eq:loja_smooth}.
\end{remark}

In principle, one could obtain a closed form expression for the last condition in \eqref{loja_strengthened_1}, by solving the associated maximization problem and obtaining an optimal $\alpha_\star$, function of $\{(x_i,f_i,g_i)\}_{i=1,2,\star}$. However, seeking a formulation that involves only linear expressions in $f_i$'s as well as in the Gram matrix of $g_i$'s,$x_i$'s, we propose instead an SDP reformulation of \eqref{loja_strengthened_1}.
\begin{proposition} \label{prop:SOS_loja_smooth}
    Consider a pair $S=\{(x_i,f_i,g_i)\}_{i=1,2}$. It holds that $S$ satisfies \eqref{loja_strengthened_1} if and only (i) $f_i\leq f_\star+\frac{\|g_i\|^2}{2 \mu }, i=1,2$, (ii) $
    f_i\geq f_\star+\frac{\|g_i\|^2}{2 L}, i=1,2$, and (iii) $\forall i,j=1,2$, there exists some
    \begin{equation} \label{def:matrix_loja}
        M^{ij}=\begin{bmatrix}
            M_{11} & M_{12}\\
            M_{12} & M_{22}
            \end{bmatrix}\succeq 0,\quad \bar M^{ij}=\begin{bmatrix}
            - A (2L+\mu) & \bar M_{12}\\
            \bar M_{12} & \bar M_{22}
        \end{bmatrix}\succeq 0,
    \end{equation}
    such that
    \begin{align*}
        M_{11}-\bar M_{11}+\frac{4\mu}{2L+\mu}\bar M_{12}&=-(B-C-(L+3\mu)A),\\
        2M_{12}-2\bar M_{12}+\frac{2\mu}{2L+\mu}\bar M_{22}&=-(L+\mu)A+2B,\\
        M_{22}-\bar M_{22}&= -B,
    \end{align*}
    where $A=\left(-f_i+f_j+\frac 12 \langle g_i+g_j,x_i-x_j\rangle +\frac{1}{4L}\|g_i-g_j\|^2-\frac{L}{4}\|x_i-x_j\|^2\right)$, $B=(L+\mu)\left(f_i-f_\star-\frac{\|g_i\|^2}{2L}\right)$, and $C=(L-\mu)\left(f_j-f_\star+\frac{\|g_j\|^2}{2L}\right)$.
\end{proposition}
The detailed proof of \cref{prop:SOS_loja_smooth} is deferred to \cref{app:SOS_loja}, and relies on Markov-Lukacs theorem, see, e.g., \cite[Theorem 1.21]{szeg1939orthogonal}.

While \eqref{def:matrix_loja} is an exact reformulation of \eqref{loja_strengthened_1}, it remains a relaxation of $\tpL$ since \eqref{loja_strengthened_1} is itself a relaxation of this one-point strengthening.

\Cref{fig:loja} illustrates, on a numerical example, the difference in datasets satisfying \eqref{eq:loja_smooth} as compared to~\eqref{loja_strengthened_1}, the relaxation of $\tpL$. Observe that in this context of simply verifying whether or not a dataset satisfies \eqref{loja_strengthened_1}, one can directly solve the maximization problem in \eqref{loja_strengthened_1}.

Alternatively, if one aims for conditions that share the same structure as \eqref{eq:loja_smooth}, one can relax \eqref{loja_strengthened_1} by imposing the last condition for $k$ specific values of $\alpha$. This could allow strengthening \eqref{eq:loja_smooth} by conducting a second one-point iteration, starting from this relaxation, see \cref{rem:complexity}.
\subsection{Convex blockwise smooth functions}\label{sec:blockwise}
This section studies the class $\F_{0,\mathbb{L}}(\R^d)$ of convex blockwise $\mathbb{L}$-smooth functions $f$, that is functions satisfying, $\forall x,y \in \R^d$ \cite{kamri2023worst},
\begin{equation} \tag{$\pblock$}
    \pblock^{xy}\leq 0
    \Leftrightarrow
    f(x)\geq f(y)+\langle \nabla f(y), x-y\rangle +\tfrac{1 }{2L_i} \|\nabla_i f(x)-\nabla_i f(y)\|^2
    ,\hspace{1cm} \ \forall i\in\llbracket n\rrbracket ,\label{eq:blockwise}
\end{equation}
where $n \in \llbracket d \rrbracket$ and $\nabla_i f$ denotes the $i^{\text{th}}$ block-component of the gradient of $f$.

It is known that this condition is a pairwise interpolation condition on $\R^d$, but no $K$-wise interpolation condition for $K\geq 3$.
\begin{example}[counterexample to sufficiency of \eqref{eq:blockwise} \cite{kamri2025worst}, Theorem 2.4]\label{ex:counterex_block}
    Let $n=d=2$, $L_1=L_2=1$, and $S=\{(x_i,f_i,g_i)\}_{i\in \llbracket 3\rrbracket}$, where $(x_1,f_1,g_1)=\left(\binom{-1}{0},\frac{1}{2}, \binom{-1}{0}\right),  \ (x_2,f_2,g_2)=\left(\binom{1}{0}, \frac{1}{2},\binom{1}{0}\right)$, and $(x_3,f_3,g_3)=\left(\binom{0}{0},0, \binom{0}{-1}\right)$. Then $S$ satisfies \eqref{eq:blockwise} while not consistent with a function in $\F_{0,\mathbb{L}}$.
\end{example}

Exploiting the one-point refinement technique, we strengthen~\eqref{eq:blockwise} on sets of cardinality $3$ as follows.
\begin{proposition} \label{prop:1pt_strengthening_blockwise}
    Let  $Q=\{(x_1,x_2,x_3)\}$, and \eqref{eq:blockwise}. Then, $\tpblock$, the terwise one-point strengthening of \eqref{eq:blockwise} with respect to $\R^d$, is satisfied only if
    \begin{align}
        \tpblock^{ijk}\leq 0 \Rightarrow
        0\geq&(1-\lambda) \bigg(-f_i+ f_j+\langle g_j, x_i-x_j\rangle+ \frac{1 }{2L_m}\|g_i^{(m)}-g_j^{(m)}\|^2\bigg)\label{eq:blockwise_strengthened}\\&+\lambda\bigg(-f_i+ f_k+\langle g_k, x_i-x_k\rangle+ \frac{1 }{2L_m}\|g_i^{(m)}-g_k^{(m)}\|^2\bigg)\nonumber\\&+\lambda(1-\lambda)\bigg(\max_l\frac{1}{2 L_l}\|g_j^{(l)}-g_k^{(l)}\|^2-\frac{1}{2L_m}\|g_j^{(m)}-g_k^{(m)}\|^2\bigg), \ \forall m\in\llbracket n\rrbracket, \ \forall \lambda \in [0,1],\nonumber
    \end{align}
    where $g_i^{(m)}$ denotes the $m^{\text{th}}$ component of $g_i$.
\end{proposition}
The detailed proof of \cref{prop:1pt_strengthening_blockwise} is deferred to \cref{app:blockwise}. It consists in (i) formulating the pairwise one-point strengthening \eqref{eq:1pt_strengthening} of $\tpblock$, (ii) dualizing the inner problem in this one-point strengthening as to relax it into a maximization problem over positive dual variables and $z\in \R^d$, (iii) explicitly maximizing this relaxation in $z$, and (iv) computing a feasible solution to this relaxation, by setting some dual variables to $0$ and obtaining a constraint involving one of the remaining variables, that is $\lambda \in [0,1]$. Observe that the $\max$ over $l$ can be easily dealt with by simply imposing \eqref{eq:blockwise_strengthened} to be true for every $l$.
\begin{remark}\label{rem:bloc}
    Observe that \cref{ex:counterex_block} is discarded by~\eqref{eq:blockwise_strengthened} with $\lambda=\frac{1}{2}$.
\end{remark}
\begin{remark}
    Setting $\lambda=0$ or $\lambda=1$ in \eqref{eq:blockwise_strengthened} allows recovering \eqref{eq:blockwise}, hence again this relaxation of $\tpblock$ is at least as strong as \eqref{eq:blockwise}. \Cref{rem:bloc} allows concluding that \eqref{eq:blockwise_strengthened} is actually strictly stronger than \eqref{eq:blockwise}.
\end{remark}
Similarly to \cref{sec:loja}, seeking a formulation that involves only linear expressions in $f_i$'s as well as in the Gram matrix of $g_i$'s,$x_i$'s, we propose a reformulation of \eqref{eq:blockwise_strengthened} as an SDP.
\begin{proposition}\label{prop:SOS_block}
    Consider a triplet $S=\{(x_i,f_i,g_i)\}_{i\in \llbracket 3\rrbracket}$. It holds that $S$ satisfies \eqref{eq:blockwise_strengthened} if and only if $\forall m\in\llbracket n\rrbracket$ and $\forall i,j,k \in \llbracket 3\rrbracket$, there exists some
    \begin{equation}
        M_m^{ijk}=\begin{bmatrix}\label{def:matrix_block}
            -A_{ijm} & M_{12} \\
            M_{12} & M_{22}
        \end{bmatrix}\succeq 0, \quad K_m^{ijk}\geq 0,
    \end{equation}
    such that
    \begin{align*}
        2M_{m,12}^{ijk}+K_m^{ijk}&=A_{ijm}-A_{ikm}-B_{jkm},\\
        M_{m,22}^{ijk}-K_m^{ijk} &=B_{jkm},
    \end{align*}
    where
    \begin{equation*}
        A_{ijm}=-f_i+ f_j+\langle g_j, x_i-x_j\rangle+ \frac{1 }{2L_m}\|g_i^{(m)}-g_j^{(m)}\|^2, \ B_{jkm}=\max_l\frac{1}{2 L_l}\|g_j^{(l)}-g_k^{(l)}\|^2-\frac{1}{2L_m}\|g_j^{(m)}-g_k^{(m)}\|^2.
    \end{equation*}
\end{proposition}
The detailed proof of \cref{prop:SOS_block} is deferred to \cref{app:SOS_block}, and relies on Markov-Lukacs theorem, see, e.g., \cite[Theorem 1.21]{szeg1939orthogonal}.

While \eqref{def:matrix_block} is an exact reformulation of \eqref{eq:blockwise_strengthened}, it remains a relaxation of $\tpblock$ since \eqref{eq:blockwise_strengthened} is itself a relaxation of this one-point strengthening.

Alternatively, optimizing \eqref{eq:blockwise_strengthened} in $\lambda$ straightforwardly allows obtaining a closed-form for this condition.
\begin{proposition}
    With the same notation as in \cref{prop:SOS_block} and considering a triplet $S=\{(x_i,f_i,g_i)\}_{i\in \llbracket 3\rrbracket}$, it holds that $S$ satisfies \eqref{eq:blockwise_strengthened} if and only if $\forall m\in\llbracket n\rrbracket$ and $\forall i,j,k \in \llbracket 3\rrbracket$,
    \begin{align}
        A_{ikm}&\leq 0 \nonumber\\
        (A_{ijm}-A_{ikm})^2+B_{jkm}^2+2B_{jkm}(A_{ijm}+A_{ikm})&\leq 0 \text{ if } \quad \frac{A_{ikm}-A_{ijm}+B_{jkm}}{2B_{jkm}} \in [0,1].\label{eq:block_quartic}
    \end{align}
\end{proposition}
\begin{proof}
    Note that \eqref{eq:blockwise_strengthened} is quadratic in $\lambda$, and possesses a maximum in
    $    \lambda_\star=\frac{A_{ikm}-A_{ijm}+B_{jkm}}{2B_{jkm}}$. When $\lambda_\star \in [0,1]$, \eqref{eq:blockwise_strengthened} is thus maximized in $\lambda_\star$ and its maximum is given by \eqref{eq:block_quartic}. Otherwise, \eqref{eq:blockwise_strengthened} is maximized in $\lambda=0$ or $\lambda=1$, leading to $A_{ikm}, \ A_{ijm}\leq 0$.
\end{proof}

\begin{figure}
    \centering
    \includegraphics[width=0.5\linewidth]{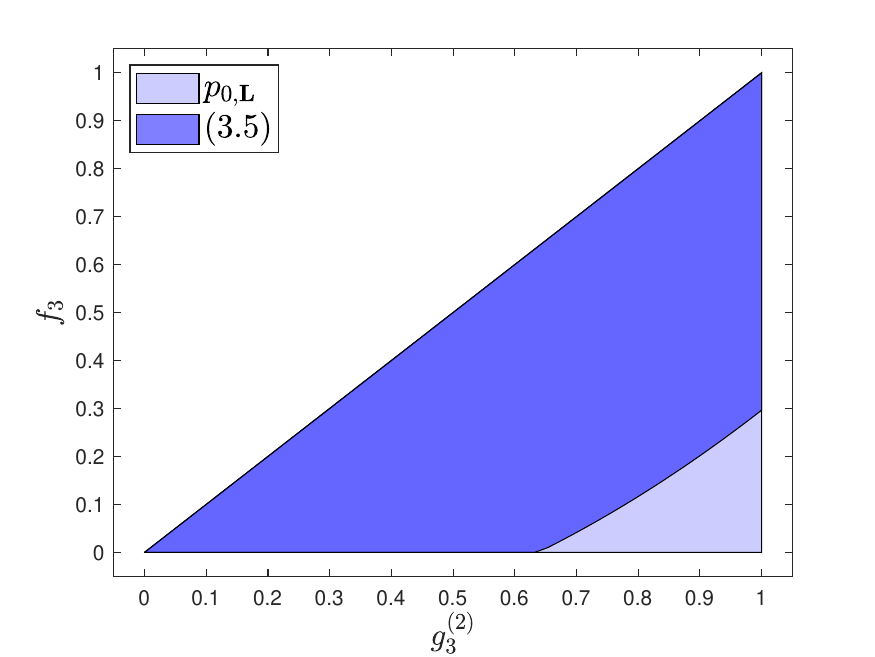}
    \caption{Allowed regions for $f_3$ as a function of $h_3$, according to \eqref{eq:blockwise_strengthened} and \eqref{eq:blockwise}, given $L_x=L_y=1$, $(x_1,f_1,g_1)=\left(\binom{-1}{0}, \frac{1}{2},\binom{-1}{0} \right), \ (x_2,f_2,g_2)=\left(\binom{1}{0},\frac{1}{2},\binom{1}{0}\right), \ x_3=\binom{0}{1}$ and $g_3=\binom{0}{g_3^{(2)}}$.}
    \label{fig:block}
\end{figure}

We illustrate in \cref{fig:block} on a numerical example the difference in datasets satisfying \eqref{eq:blockwise} as compared to~\eqref{eq:blockwise_strengthened}, the relaxation of $\tpblock$.

Similarly to \cref{sec:loja}, one can obtain conditions that share exactly the same structure as \eqref{eq:blockwise} by imposing \eqref{eq:blockwise_strengthened} for $k$ specific values of $\lambda$, which could provide a starting point for a second iteration of the refinement procedure, see \cref{rem:complexity}.
\subsection{Strongly monotone Lipschitz operators}\label{sec:operators}
This section studies the class $\Op(\R^d)$ of maximally monotone Lipschitz operators, that is operators $T$ satisfying, for some $0\leq\mu\leq L$, and $\forall x,y \in \R^d$ \cite{bauschke2017convex,ryu2020operator},
\begin{equation}\tag{$\pop$}
    \pop^{xy}\leq 0 \Leftrightarrow \begin{cases}
        - \langle T(x)-T(y),x-y\rangle+\mu\|x-y\|^2& \leq 0\\
        \quad\quad\quad\ \|T(x)-T(y)\|- L \|x-y\|&\leq 0.
    \end{cases} \hspace{2cm} \label{eq:operators}
\end{equation}
To prevent any confusion with properties related to function classes, we use the notation $q$ instead of $p$ to denote a property.

It is known \cite[Proposition 4]{ryu2020operator} that \eqref{eq:operators} is a pairwise interpolation condition on $\R^d$, and a global interpolation condition when $d=1$. However, when $K\geq 3$ or $d\geq 2$, \eqref{eq:operators} is no longer a $K$-wise interpolation condition and is only necessary.
\begin{example}[counterexample to sufficiency of \eqref{eq:operators} \cite{ryu2020operator}]
    Let $d=2$, $\mu=0$, $L=1$, and $S=\{(x_i,t_i)\}_{i\in \llbracket 3\rrbracket}$, where $(x_1,t_1)=\left(\binom{0}{0},\binom{0}{0}\right), \ (x_2,t_2)=\left(\binom{1}{0},\binom{0}{0}\right), \ (x_3,t_3)=\left(\binom{1/2}{0},\binom{0}{-1/2}\right)$. Then $S$ satisfies \eqref{eq:operators} while not consistent with a function in $\Op$.\label{ex:op}
\end{example}

The refinement procedure described in \cref{sec:formal_extension} allows strengthening \eqref{eq:operators}, starting with a set of cardinality $3$.
\begin{proposition} \label{prop:1pt_strengthening_operators}
    Consider $Q=\{(x_1,x_2,x_3)\}$, and \eqref{eq:operators}.
    Then, $\tpop$, the terwise one-point strengthening of \eqref{eq:operators} with respect to $\R^d$, is satisfied only if
    \begin{equation}
        \tpop^{ijk}\leq 0 \Rightarrow\left\{
            \begin{aligned}
                &\forall  \lambda_k, \ \mu_i, \ \mu_j \geq 0:\\
                &0\geq
                \lambda_k(\mu_i+\mu_j)(A_{ik}\mu_i+A_{jk}\mu_j)-\lambda_k \mu_i\mu_j A_{ij}+\lambda_k^2 (B_{ik} \mu_i +B_{jk}\mu_j)
                \\\MoveEqLeft[-1]+B_{ij} \mu_i\mu_j(\mu_i+\mu_j-2\mu\lambda_k)\\
                &0\geq \lambda_k(\mu_i+\mu_j)(B_{ik}\mu_i+B_{jk}\mu_j)-\lambda_k \mu_i\mu_j B_{ij}+\lambda_k^2 (A_{ik} \mu_i+A_{jk}\mu_j)
                \\\MoveEqLeft[-1] +A_{ij} \mu_i\mu_j(\mu_i+\mu_j-2\mu \lambda_k),
        \end{aligned}\right.
        \label{eq:final_p_op}
    \end{equation}
    where $A_{ij}=\|t_i-t_j\|^2-L^2\|x_i-x_j\|^2$ and $B_{ij}=-2L\langle t_i-t_j,x_i-x_j\rangle+2L\mu\|x_i-x_j\|^2$.
\end{proposition}
The detailed proof of \cref{prop:1pt_strengthening_operators} is deferred to \cref{app:op}.  It consists in (i) formulating the pairwise one-point strengthening \eqref{eq:1pt_strengthening} of $\tpop$, (ii) dualizing the inner problem in this one-point strengthening as to relax it into a maximization problem over positive dual variables and $z\in \R^d$, (iii) explicitly maximizing this relaxation in $z$, and (iv) computing a feasible solution to this relaxation, by setting some dual variables to $0$ and obtaining a constraint involving the remaining variables, i.e., $\lambda_k, \mu_i$ and $\mu_j$.

\begin{remark}\label{rem:op}
    Observe that \cref{ex:op} does not satisfy \eqref{eq:final_p_op}, and separates $\F_{\pop}(Q)$ from $\F_{\tpop}(Q)$.
\end{remark}
\begin{remark}
    Setting $\lambda_k=0$ in \eqref{eq:final_p_op} allows recovering \eqref{eq:operators}, hence again this relaxation of $\tpop$ is at least as strong as \eqref{eq:operators}. \Cref{rem:op} allows concluding \eqref{eq:final_p_op} is actually strictly stronger than \eqref{eq:operators}.
\end{remark}

Similarly to \cref{sec:loja,sec:blockwise}, we propose a relaxation of \eqref{eq:final_p_op} as an SDP.
\begin{proposition}\label{prop:SOS_op}
    With the same notation as in \cref{prop:1pt_strengthening_operators} and considering a triplet $S=\{(x_i,t_i)\}_{i\in \llbracket 3\rrbracket}$, it holds that $S$ satisfies \eqref{eq:final_p_op} if ~$\forall i,j,k \in \llbracket 3\rrbracket$, there exists some $M_{A,B}^{ijk}\succeq0$, $M_{B,A}^{ijk}\succeq0$, where
    \begin{equation}
        M_{\alpha,\beta}^{ijk}=\begin{bmatrix}\label{def:matrix_op}
            -\beta_{ij} & 0 & 0 & M_{14} &M_{15}& M_{16}&M_{17}\\
            0 & -\alpha_{jk} & -M_{15} & 0 & -M_{14}&M_{26}&M_{27}\\
            0 & -M_{15} & -\beta_{ij}  & M_{34}& -M_{16}&0&M_{37}\\
            M_{14} & 0 & M_{34} &-\beta_{jk}&-M_{27}&M_{46}&0\\
            M_{15}&-M_{14} &-M_{16} & -M_{27}&M_{55}&-M_{37}&-M_{46}\\
            M_{16}&M_{26} &0 &M_{46} & -M_{37}&-\alpha_{ik}&0\\
            M_{17}&M_{27} &M_{37} &0 &-M_{46} & 0& -\beta_{ik}
        \end{bmatrix},
    \end{equation}
    and
    \begin{equation*}
        M_{55}+2M_{17}+2M_{26}+2M_{34}=\alpha_{ij}+2\mu\beta_{ij}-\alpha_{jk}-\alpha_{ik}.
    \end{equation*}
\end{proposition}
The detailed proof of \cref{prop:SOS_op} is deferred to \cref{app:SOS_op}.

\begin{figure}
    \centering
    \includegraphics[width=0.5\linewidth]{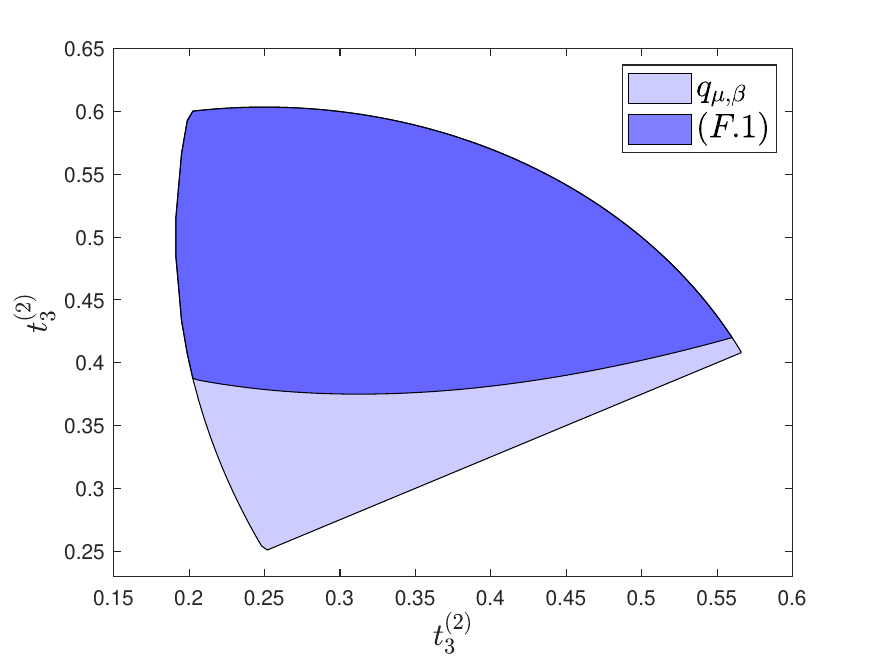}
    \caption{Allowed regions for $t_3$, according to \eqref{eq:operators} or \eqref{eq:final_p_op}, given $\mu=0$, $L=1$, $(x_1,t_1)=\left(\binom{0}{0},\binom{0}{0}\right)$, $(x_2,t_2)=\left(\binom{1}{0},\binom{0}{-1/2}\right)$, and $x_3=\binom{1/2}{1/2}$.}
    \label{fig:operators}
\end{figure}

We illustrate in \cref{fig:operators} on a numerical example the difference in datasets satisfying~\eqref{eq:operators} as compared to~\eqref{eq:final_p_op}, the relaxation of $\tilde{q}_{0,L}$.

\section{Concluding remarks}
This work presents a one-point strengthening method, that can be seen as a principled and constructive approach to interpolation/extension conditions. The approach starts from any characterization of a class of functions/operators and constructs a refined characterization of this class. Repeating the procedure iteratively allows obtaining extension conditions as a final product. In particular, it allows for strengthening given sets of conditions based only on their algebraic properties, leading to improved conditions that would have been challenging to guess based solely on the class' analytical properties. In addition, starting from linear conditions in function values and Gram representation of points and subgradients (or operator values) yields refined conditions that are SDP-representable and can be used as such when imposed on datasets that are variable to some optimization problem.

The approach naturally showcases some limitations and raises open questions. First, our main technical tools for computing one-point strengthening consist of dualizing the inner problem in~\eqref{eq:1pt_strengthening}. Hence, the approach fails as soon as the initial property under consideration is not convex in $F_z$. Similarly, the approach is naturally limited to properties for which this dualization yields a zero duality gap. In other cases, the approach might attain only a very loose strengthening. Furthermore, we generally have to work with closed-form solutions to this inner problem, and we are thereby quickly limited to having to deal with reasonably simple expressions. Indeed, we quickly arrive at formulations that turn out to be hardly manipulable, both for computing their one-point strengthening and for proving their stability.

The many remaining open challenges thus include the questions of tackling conditions whose one-point strengthenings are not convex in $F_z$. Another venue for improvement is that of specializing the approach to more specific classes of problems (e.g., smooth functions that are defined on specific sets instead of smooth functions with general convex domains; see \cref{sec:extension_subset}). Yet another possible improvement could be to search for one-point strengthenings without relaxations, for instance, by relying on quantifier elimination techniques (see, e.g.,\cite{algorithmsrealalgebraic}) -- that do not scale well in the number of variables and parameters in general, but that could handle our problems with very few of them.
\paragraph{Codes.} Codes for reproducing and verifying technical parts of this work can be found at
\begin{center}
    \url{https://github.com/AnneRubbens/Constraints_Strengthening}
\end{center}
As this work is strongly motivated by algorithm analyses -- it was recently shown that interpolation/extension conditions are key for systematically analyzing their performance; see, e.g.,~\cite{taylor2017smooth} or recent presentations~\cite{goujaud2023fundamental,rubbens2023interpolation} --  we provide a few examples of such applications using the conditions obtained in \cref{sec:applications}. These include analysis of the gradient method on smooth function satisfying a \L{}ojasiewisz condition~\cite{bolte2017error,abbaszadehpeivasti2022conditionslinearconvergencegradient,karimi2016linear} (see \cref{sec:loja}), of the random or cyclic coordinate descent
on convex blockwise smooth functions \cite{abbaszadehpeivasti2022convergence,beck,kamri2023worst,nesterov2012efficiency} (see \cref{sec:blockwise}), and of the optimistic gradient and extragradient on Lipschitz monotone operators \cite{gorbunov2022last,gorbunov2022extragradient,korpelevich1976extragradient} (see \cref{sec:operators}). In addition, the SDP formulations \eqref{def:matrix_loja}, \eqref{def:matrix_block} and \eqref{def:matrix_op} of the improved descriptions of (i) smooth functions satisfying a \L{}ojsiewicz condition, (ii) blockwise smooth convex functions, and (iii) (strongly) monotone Lipschitz operators, have been added to the PEPit Python library \href{https://pepit.readthedocs.io/en/latest/api/functions_and_operators.html}{PEPit documentation}~\cite{goujaud2024pepit}, which allows solving PEPs on these classes straightforwardly.

\appendix

\section{Showing stability of a condition on an example: Lipschitz operators}\label{app:op_Lip}
Consider the class $\mathcal{Q}_{-L,L}$ of Lipschitz operators $T$, satisfying, $\forall x, y \in \R^d$,
\begin{equation} \tag{$q_{-L,L}$}\label{eq:Lipschitz_op}
    q_{-L,L}^{xy}\leq 0\Leftrightarrow \|T(x)-T(y)\|\leq L \|x-y\|.
\end{equation}
It is known \cite{kirszbraun1934zusammenziehende,valentine1943extension} that \eqref{eq:Lipschitz_op} is a global interpolation constraint for $\mathcal{Q}_{-L,L}$.

This section shows how to retrieve this result via the refinement procedure described in \cref{sec:formal_extension}, by analyzing different formulations of the one-point strengthening of \eqref{eq:Lipschitz_op}.

Consider $Q=\{(x_i)\}_{i\in \llbracket N\rrbracket}$. The $N$-wise strengthening of \eqref{eq:Lipschitz_op} with respect to $\R^d$ is given by
\begin{equation*}
    \tilde{q}_{-L,L}= \max_{z\in \R^d}\min_{\substack{\tau \in \R, \\ t_z\in \R^d}} \tau
    \text{ s.t. }\|t_z-t_i\|^2-L^2\|z-x_i\|^2\leq \tau, \ i\in \llbracket N\rrbracket, \hspace{2cm}  (\lambda_i)
\end{equation*}
where $\lambda_i$ are the associated coefficients. Observe that this strengthening is convex in $t_z$ and $\tau$, and that $(t_z,\tau)=(0,\max_i \|t_i\|^2+1)$ satisfy Slater's condition, hence there is no duality gap for this inner problem. Dualization of the inner problem yields
\begin{equation*}
    \tilde{q}_{-L,L}=  \max_{\substack{\lambda_i \geq 0, \\ \sum_i\lambda_i=1}} \  \sum_{i} \lambda_i (\|t_z^\star-t_i\|^2-L^2\|z-x_i\|^2) \text{ s.t. } t_z^\star=\sum_i\lambda_i t_i \text{ and } z=\sum_i \lambda_i x_i,
\end{equation*}
where $t_z^\star$ minimizes the inner problem in $\tilde{q}_{-L,L}$. This condition is homogeneous of degree $2$ in $\lambda_i$. Hence, for any $K\geq 0$, it must hold with $\tilde \lambda_i=K\lambda_i$ replacing $\lambda_i$, which allows removing the condition $\sum_i \lambda_i =1$.

Equivalently, it thus holds
\begin{equation}\label{eq:op_ex}
    \tilde{q}_{-L,L}= \max_{\lambda_i \geq 0} \sum_{i,j\in \llbracket N\rrbracket, \ j>i} \lambda_j \lambda_i (\|t_i-t_j\|^2-L^2\|x_i-x_j\|^2).
\end{equation}

Necessarily, for $\tilde{q}_{-L,L}$ to be negative, it must hold $\|t_i-t_j\|^2-L^2\|x_i-x_j\|^2\leq 0$, $i,j \in \llbracket N\rrbracket$. In addition, if this condition is satisfied, then $\tilde{q}_{-L,L}\leq0$. Hence, $\tilde{q}_{-L,L}=q_{-L,L}$.
\section{Smooth functions satisfying a quadratic \L{}ojasiewicz condition}
\subsection{Derivation of the one-point strengthening: proof of Proposition \ref{prop:1pt_strengthening_loja_smooth}} \label{app:loja_smooth}
\begin{proof}
    Without loss of generality, we let $f_\star=0,x_\star=0$. The pairwise one-point strengthening of \eqref{eq:loja_smooth} on $\R^d$ is given by
    \begin{equation*}
        \begin{aligned}
            \tpL= \max_{z\in \R^d}\min_{\substack{\tau \in \R, \\ g_z\in \R^d,\ f_z\in \R}} & \tau && \\
            \text{ s.t. } & f_z-f_i+\frac{1}{2}\langle g_z+ g_i,x_i-z\rangle-\frac{L}{4}\|z-x_i\|^2+\frac{1}{4L}\|g_z-g_i\|^2\leq \tau, \quad i=1,2,\star&&    (\mu_i)\\
            &f_i-f_z+\frac{1}{2}\langle g_i+ g_z,z-x_i\rangle-\frac{L}{4}\|z-x_i\|^2+\frac{1}{4L}\|g_z-g_i\|^2\leq \tau, \quad i=1,2,\star &&(\lambda_i) \\
            &f_z-\frac{1}{2\mu} \|g_z\|^2\leq  \tau && (\alpha_1)\\
            &-f_z\leq \tau,&&(\alpha_2)\\
        \end{aligned}
    \end{equation*}
    where $\alpha_1, \ \alpha_2$, and $\lambda_i, \mu_i, \ i=1,2,\star$ are the dual coefficients. Partial dualization of the inner problem over $\tau$ and $f_z$, and minimization over $g_z$ yields
    \begin{equation*}
        \begin{aligned}
            \tpL \geq \max_{\substack{z\in \R^d,\\ \lambda_i\geq 0,\ \mu_i\geq 0, \\ \alpha_1, \ \alpha_2\geq 0 }} &\sum_{i\in \mathcal{I}} \big(\lambda_i (f_i+\frac{1}{2}\langle g_i+ g_z^\star,z-x_i\rangle-\frac{L}{4}\|z-x_i\|^2+\frac{1}{4L}\|g_z^\star-g_i\|^2)  \\&\quad+ \mu_i (-f_i+\frac{1}{2}\langle g_z^\star+ g_i,x_i-z\rangle-\frac{L}{4}\|z-x_i\|^2+\frac{1}{4L}\|g_z^\star-g_i\|^2)\big)\nonumber\\&\quad -\alpha_1 \frac{\|g_z^\star\|^2}{2\mu} \\
            \text{ s.t. }
            \frac{1}{2}&=\sum_{i\in \mathcal{I}} \lambda_i+\alpha_2=\sum_{i\in \mathcal{I}} \mu_i+\alpha_1,\\
            0& =\sum_{i\in \mathcal{I}}\left((\mu_i+\lambda_i) \frac{1}{2L} (g_z^\star-g_i)+\frac{\lambda_i-\mu_i}{2} (z-x_i)\right) -\alpha_1 \frac{g_z^\star}{\mu},\\
            0&<\sum_{i\in \mathcal{I}}\frac{\mu_i+\lambda_i }{2L}  - \frac{\alpha_1}{\mu},
        \end{aligned}
    \end{equation*}
    where $\mathcal{I}=\{1,2,\star\}$, $g_z^\star$ minimizes the inner problem in $\tpL$, and the last condition ensures that $g_z^\star$ is a minimizer and no maximizer.

    In particular, setting all dual coefficients to $0$ except for $\alpha_2=\mu_i=\frac{1}{2}$, for a given $i\in \mathcal{I}$, $z=x_i-\frac{g_i}{L}$ and $g_z^\star=0$ leads to a feasible solution of $-f_i+\frac{\|g_i\|^2}{2L}$, hence \eqref{loja_strengthened_1} cannot be satisfied unless this quantity is negative.

    In addition, setting all dual coefficients to $0$ except for
    $\alpha_1\leq\frac{\mu}{2L+\mu}, \quad \lambda_j=\frac{1}{2}, \quad \mu_i=\frac{1}{2}-\alpha_1$,
    for any $i,j\in \mathcal{I}$, and maximizing over $z$, leads to a feasible solution \begin{align*}
        \tpL\geq 0 &\Leftrightarrow 0\geq-f_i+ f_j+\frac{1}{2}\langle g_i+g_j, x_i-x_j\rangle+\frac{1}{4L}\|g_i-g_j\|^2-\frac{L}{4} \|x_i-x_j\|^2\\ &+\frac{\alpha_1}{(1-2\alpha_1)(\mu-(L+\mu)\alpha_1)} \left((1-2\alpha_1)^2(L+\mu)(f_i-f_\star-\frac{\|g_i\|^2}{2L})-(L-\mu)(f_j-f_\star+\frac{\|g_j\|^2}{2L})\right).
    \end{align*}
    Hence, $\tpL$ cannot be satisfied unless this quantity is negative for all $0\leq\alpha_1\leq\frac{\mu}{2L+\mu}$, and in particular for its maximal value with respect to $\alpha_1$. Setting $\alpha=2\alpha_1$, it holds that $\tpL$ is at least as strong as \eqref{loja_strengthened_1}.
\end{proof}
\subsection{Derivation of a SOS reformulation: proof of Proposition \ref{prop:SOS_loja_smooth}}\label{app:SOS_loja}
\begin{proof}
    The first conditions are strictly equivalent to the first ones in \eqref{loja_strengthened_1}. The last inequality in \eqref{loja_strengthened_1} can be equivalently written as \begin{equation}
        0\geq A+ \frac{\alpha}{(1-\alpha)(2\mu-(L+\mu)\alpha)} \left((1-\alpha)^2B-C)\right), \ \forall \alpha \in [0,\frac{2\mu}{2L+\mu}],\label{eq1}
    \end{equation}
    where
    \begin{align*}
        A&=\left(-f_i+f_j+\frac 12 \langle g_i+g_j,x_i-x_j\rangle +\frac{1}{4L}\|g_i-g_j\|^2-\frac{L}{4}\|x_i-x_j\|^2\right)\\
        B&=(L+\mu)\left(f_i-f_\star-\frac{\|g_i\|^2}{2L}\right)\\
        C&=(L-\mu)\left(f_j-f_\star+\frac{\|g_j\|^2}{2L}\right).
    \end{align*}
    In the interval $[0,\frac{2\mu}{2L+\mu}]$, $\alpha \leq \frac{2\mu}{L+\mu} \leq 1$, hence \eqref{eq1} is equivalent to
    \begin{align*}
        &0\geq (1-\alpha)(2\mu-(L+\mu)\alpha)A+ (1-\alpha)^2\alpha B-\alpha C, \ \forall \alpha \in [0,\frac{2\mu}{2L+\mu}]\\
        \Leftrightarrow\ &0\geq 2\mu A+\alpha (B-C-(L+3\mu)A)+\alpha^2((L+\mu)A-2B)+\alpha^3B, \ \forall \alpha \in [0,\frac{2\mu}{2L+\mu}]\\
        \Leftrightarrow\ &0\geq a_0+a_1\alpha +a_2\alpha^2+a_3\alpha^3:=P(\alpha), \ \forall \alpha \in [0,\frac{2\mu}{2L+\mu}],
    \end{align*}
    where $a_0=2\mu A$, $a_1=B-C-(L+3\mu)A$, $a_2=(L+\mu)A-2B$ and $a_3=B$.

    By Markov-Lukacs theorem, see, e.g., \cite[Theorem 1.21]{szeg1939orthogonal},
    \begin{equation*}
        P(\alpha)\leq 0, \forall \alpha \in [0,\frac{2\mu}{2L+\mu}] \Leftrightarrow \exists M, \bar M \succeq 0 \text{ s.t. } P(\alpha)=-\alpha \begin{bmatrix}
            1&\alpha
            \end{bmatrix} M \begin{bmatrix}
            1\\ \alpha
            \end{bmatrix}- (\frac{2\mu}{2L+\mu}-\alpha) \begin{bmatrix}
            1&\alpha
            \end{bmatrix} \bar M \begin{bmatrix}
            1\\ \alpha
        \end{bmatrix}.
    \end{equation*}
    This condition can be stated linearly in terms of the coefficients of $M,\bar M$ and those of $P(\cdot)$. Denoting the $i,j^{\text{th}}$ entrance of $M$ by $M_{ij}$, $i,j\in [ 2]$, the list of constraints is given by
    \begin{equation*}
        \begin{aligned}[b]
            \frac{2\mu}{2L+\mu}\bar M_{11}&=-a_0,&
            M_{11}-\bar M_{11}+\frac{4\mu}{2L+\mu}\bar M_{12}&=-a_1\\
            2M_{12}- 2 \bar M_{12}+\frac{2\mu}{2L+\mu}\bar M_{22} &=-a_2, &
            M_{22}-\bar M_{22}&=-a_3.
        \end{aligned}
        \qedhere
    \end{equation*}
\end{proof}
\section{Convex blockwise smooth functions}
\subsection{Derivation of the one-point strengthening: proof of Proposition \ref{prop:1pt_strengthening_blockwise}}\label{app:blockwise}
\begin{proof}
    The terwise one-point strengthening of \eqref{eq:blockwise} on $\R^d$ is given by
    \begin{align*}
        \tpblock= \max_{z\in \R^d}\min_{\substack{\tau \in \R, \\ g_z\in \R^d,\ f_z\in \R}} &\tau && \\
        \text{ s.t. }   &f_z-f_i+\langle g_z,x_i-z\rangle+\frac{1}{2L_m}\|g_z^{(m)}-g_i^{(m)}\|^2\leq  \tau,\quad i\in \llbracket 3\rrbracket, \ m\in\llbracket n\rrbracket&&
        (\mu_{i,m})\nonumber\\
        &f_i-f_z+\langle g_i,z-x_i\rangle+\frac{1}{2L_m}\|g_z^{(m)}-g_i^{(m)}\|^2\leq  \tau, \quad i\in \llbracket 3\rrbracket, \ m\in\llbracket n\rrbracket &&(\lambda_{i,m}), \nonumber
    \end{align*}
    where $\lambda_{i,m}, \mu_{i,m}, \ i\in \llbracket 3\rrbracket, \ m \in \llbracket n \rrbracket$, are the associated dual coefficients. Partial dualization of the inner problem over $\tau$ and $f_z$, and minimization over $g_z$ yields
    \begin{align*}
        \tpblock= \max_{\substack{z\in \R^d, \\ \mu_{i,m}\geq 0,\\ \lambda_{i,m}\geq 0} }&\sum_i\sum_m \mu_{i,m}(-f_i+\langle g_z^\star,x_i-z\rangle\\&+\frac{1}{2L_m}\|g_z^{\star(m)}-g_i^{(m)}\|^2)+\lambda_{i,m}(f_i+\langle g_i,z-x_i\rangle+\frac{1}{2L_m}\|g_z^{\star(m)}-g_i^{(m)}\|^2)\nonumber\\
        \text{s.t. } &\sum_i\sum_m \mu_{i,m}=\sum_i\sum_m \lambda_{i,m}=\frac{1}{2}, \text{ and } \forall m\in \llbracket n \rrbracket,\nonumber\\
        &\sum_i \sum_m\left(\mu_{i,m}(x_i^{(m)}-z^{(m)})+\mu_{i,m}\frac{g_z^{\star(m)}-g_i^{(m)}}{L_m}+\lambda_{i,m}\frac{g_z^{\star(m)}-g_i^{(m)}}{L_m}\right)=0\nonumber,
    \end{align*}
    where $g_z^\star$ minimizes the inner problem in $\tpblock$.
    In particular, setting all dual coefficients to $0$ except for $\mu_{i,m}=\frac{1}{2}$, $\lambda_{k,l}$ and $\lambda_{j,l}=\frac{1}{2}-\lambda_{k,l}$ for some $i,j \in \llbracket 3\rrbracket$ and $m,l\in \llbracket n\rrbracket$ fixed yields
    \begin{equation*}
        \begin{aligned}
            0&\geq(\frac{1}{2}-\lambda_{k,l}) (-f_i+ f_j+\langle g_j, x_i-x_j\rangle+ \frac{1 }{2L_m}\|g_i^{(m)}-g_j^{(m)}\|^2)\\\MoveEqLeft[-1] +\lambda_{k,l}(-f_i+ f_k+\langle g_k, x_i-x_k\rangle+ \frac{1 }{2L_m}\|g_i^{(m)}-g_k^{(m)}\|^2)\nonumber\\\MoveEqLeft[-1]+\lambda_{k,l}(1-2\lambda_{k,l})\left(\max_l\frac{1}{2 L_l}\|g_j^{(l)}-g_k^{(l)}\|^2-\frac{1}{2L_m}\|g_j^{(m)}-g_k^{(m)}\|^2\right), \ \forall m\in\llbracket n\rrbracket, \ \forall \lambda_{k,l} \in [0,\frac{1}{2}].\nonumber
        \end{aligned}
    \end{equation*}
    Setting $\lambda=2\lambda_{k,l}$ yields \eqref{eq:blockwise_strengthened}, hence $\tpblock$ is at least as strong as \eqref{eq:blockwise_strengthened}.
\end{proof}
\subsection{Derivation of a SOS reformulation: proof of Proposition \ref{prop:SOS_block}}
\label{app:SOS_block}
\begin{proof}
    Condition \eqref{eq:blockwise_strengthened} can be equivalently expressed as\begin{equation*}
        0\geq(1-\lambda) A_{ijm}+\lambda A_{ikm}+\lambda(1-\lambda)B_{jkm}:=P(\lambda), \ \forall m\in\llbracket n\rrbracket, \ \forall \lambda \in [0,1],
    \end{equation*}
    where
    \[
        A_{ijm}=-f_i+ f_j+\langle g_j, x_i-x_j\rangle+ \frac{1 }{2L_m}\|g_i^{(m)}-g_j^{(m)}\|^2\leq 0
    \]
    and
    \[
        B_{jkm}=\max_l\frac{1}{2 L_l}\|g_j^{(l)}-g_k^{(l)}\|^2-\frac{1}{2L_m}\|g_j^{(m)}-g_k^{(m)}\|^2\geq 0.
    \]
    By Markov-Lukasc theorem \cite[Theorem 1.21]{szeg1939orthogonal}, 
    \begin{equation*}
        P(\lambda)\leq 0, \forall \lambda \in [0,1] \Leftrightarrow \exists M\succeq 0, K\geq 0 \text{ s.t. } P(\lambda)=-\begin{bmatrix}
            1&\lambda
            \end{bmatrix} M \begin{bmatrix}
            1\\ \lambda
        \end{bmatrix}-\lambda(1-\lambda) K.
    \end{equation*}
    This condition can be stated linearly in terms of the coefficients of $M$, of $K$ and of the coefficients of $P(\cdot)$. Denoting the $i,j^{\text{th}}$ entrance of $M$ by $M_{ij}$, $i,j\in [2]$, the list of constraints is given by
    \begin{equation*}
        \begin{aligned}
            M_{11}&=-A_{ijm}, \quad
            2M_{12}+K=A_{ijm}-A_{ikm}-B_{jkm}, \quad
            M_{22}-K =B_{jkm}.
        \end{aligned}
        \qedhere
    \end{equation*}
\end{proof}
\section{Strongly monotone Lipschitz operators}
\subsection{Derivation of the one-point strengthening: proof of Proposition \ref{prop:1pt_strengthening_operators}}\label{app:op}
\begin{proof}
    The terwise one-point strengthening of \eqref{eq:operators} on $\R^d$ is given by
    \begin{align*}
        \tpop= \max_{z\in \R^d}\min_{\substack{\tau \in \R, \\ t_z\in \R^d}} \tau
        \text{ s.t. }& \frac{\|t_z-t_i\|^2}{2L}-\frac{L}{2}\|z-x_i\|^2\leq \tau, \ i\in \llbracket 3\rrbracket  && (\lambda_i) \nonumber \\
        & -\langle t_z-t_i,z-x_i\rangle+\mu\|z-x_i\|^2\leq \tau, \   i\in \llbracket 3\rrbracket&& (\mu_i) \nonumber
    \end{align*}
    where $\lambda_i,\mu_i$ are the associated dual coefficients. Partial dualization of the inner problem over $\tau$, and minimization over $t_z$ yields
    \begin{align*}
        \tpop\geq \frac{1}{2L}\max_{\substack{ \lambda_i\geq 0, \\ \mu_i \geq 0}}&\sum_{i}\big  (-\lambda_i\left(\|t_z^\star-t_i\|^2-L^2\|z-x_i\|^2\right) +\mu_i 2L(-\langle t_z^\star-t_i,z-x_i\rangle+\mu\|z-x_i\|^2) \big)\\
        \text{ s.t. }\ \  &\sum_i \frac{1}{L}\lambda_i (t_z^\star-t_i)-\mu_i (z-x_i)=0, \sum_i (\lambda_i+\mu_i) =1,
    \end{align*}
    where $t_z^\star$ minimizes the inner problem in $\tpop$. Maximizing over $z$, multiplying by $2L\big ( (\sum_i \lambda_i)^2+ (\sum_i \mu_i)^2\big)$ and rearranging the terms yields $\tpop\geq 0$ if and only if
    \begin{align*}
        0\geq\max_{ \substack{\lambda_i\geq 0, \ \mu_i \geq 0, \\ \sum_i (\lambda_i+\mu_i) =1}} \  \sum_{i,j \in \llbracket 3\rrbracket:\ j>i} \bigg(
            &A_{ij} \big(\lambda_i\lambda_j\sum_k \lambda_k+\lambda_i \mu_j (\sum_k \mu_k)+\lambda_j \mu_i (\sum_k \mu_k)-(\sum_k \lambda_k)\mu_i\mu_j\\
            \MoveEqLeft[-17.5]-2\mu(\sum_k\mu_k)\lambda_i\lambda_j\big)
            \\+& B_{ij}\big(\mu_i\mu_j\sum_k \mu_k+\mu_i \lambda_j (\sum_k \lambda_k)+\mu_j \lambda_i (\sum_k \lambda_k)-(\sum_k \mu_k)\lambda_i\lambda_j\\
        \MoveEqLeft[-17.5]-2\mu(\sum_k\lambda_k)\mu_i\mu_j\big)\bigg),
    \end{align*}
    where $A_{ij}=\|t_i-t_j\|^2-L^2\|x_i-x_j\|^2$ and $B_{ij}=-2L\langle t_i-t_j,x_i-x_j\rangle+2L\mu\|x_i-x_j\|^2$.

    This condition is homogeneous of degree $3$ in $\lambda_k, \ \mu_i,\ \mu_j$. Hence, for any $K\geq 0$, it must hold with $\tilde \mu_i=K\mu_i$, $\tilde \lambda_i=K\lambda_i$ replacing $\mu_i$, $\lambda_i$ respectively, which allows removing the condition $\sum_i (\lambda_i+\mu_i) =1$.

    Setting only $\lambda_k,\ \mu_i$ and $\mu_j$ non-zero yields the first condition in \eqref{eq:final_p_op}, and setting only $\mu_k,\ \lambda_i$ and $\lambda_j$ non-zero yields the second condition in \eqref{eq:final_p_op}, hence $\tpop$ is at least as strong as \eqref{eq:final_p_op}.
\end{proof}

\subsection{Derivation of a SOS reformulation: proof of Proposition \ref{prop:SOS_op}}
\label{app:SOS_op}
\begin{proof}
    Take the first condition in \eqref{eq:final_p_op}, and set $\lambda_k=\gamma^2$, $\mu_i=\theta^2$ and $\mu_j=\alpha^2$. The condition becomes
    \begin{equation*}
        0\geq
        \gamma^2(\alpha^2+\theta^2)(A_{ik}\theta^2+A_{jk}\alpha^2)-\gamma^2\theta^2\alpha^2 A_{ij}+\gamma^4 (B_{ik} \theta^2+B_{jk}\alpha^2)+B_{ij} \theta^2\alpha^2(\theta^2+\alpha^2-2\mu\gamma^2)=P(\alpha,\theta,\gamma).
    \end{equation*}
    $P(\alpha,\theta,\gamma)$ is an homogeneous polynomial of $3$ variables and degree $6$, which must be negative for all $\theta,\alpha,\gamma\in\R$. This negativity condition can be relaxed by requiring $P$ to be Sum of Squares, that is, by requiring the existence of some $M\preceq 0$ such that
    \[
        P(\alpha, \theta, \gamma)=\begin{pmatrix}
            \alpha^2\theta\\\alpha^2\gamma\\ \alpha \theta^2\\\alpha \gamma^2\\\alpha\theta \gamma \\ \theta^2\gamma\\\theta\gamma^2
            \end{pmatrix}^\top M \begin{pmatrix}
            \alpha^2\theta\\\alpha^2\gamma\\ \alpha \theta^2\\\alpha \gamma^2\\\alpha\theta \gamma \\ \theta^2\gamma\\\theta\gamma^2
        \end{pmatrix},
    \]
    which can be stated linearly in terms of the coefficients of $M$ and those of $P(\cdot)$. Denoting the $i,j^{\text{th}}$ entrance of $M$ by $M_{i,j}$, $i,j\in \llbracket 7\rrbracket$, the list of constraints is given by
    \begin{equation*}
        \begin{aligned}
            M_{1,1}=M_{3,3}=B_{ij}\\
            M_{2,2}=A_{jk}\\
            M_{4,4}=B_{jk}\\
            M_{6,6}=A_{ik}\\
            M_{7,7}=B_{ik}\\
            M_{5,5}+2M_{1,7}+2M_{2,6}+2M_{3,4}=-A_{ij}+A_{ik}+A_{jk}-2\mu B_{ij}\\
            M_{1,2}=M_{1,3}=M_{2,4}=M_{3,6}=M_{4,7}=M_{6,7}=0\\
            M_{1,5}+M_{2,3}=M_{1,4}+M_{2,5}=M_{1,6}+M_{3,5}=M_{2,7}+M_{4,5}=M_{3,7}+M_{5,6}=M_{4,6}+M_{5,7}=0.
        \end{aligned}
    \end{equation*}

    Translating the second condition in \eqref{eq:final_p_op} follows exactly the same argument, with $A_{(\cdot)}$ and $B_{(\cdot)}$ inverted.
\end{proof}
\section{Strongly monotone cocoercive operators}\label{sec:operators2}

This section studies the class $\Opi(\R^d)$ of strongly monotone cocoercive operators, that is operators $T$ satisfying, $\forall x,y \in \R^d$ \cite{bauschke2017convex,ryu2020operator},
\begin{equation}\tag{$\popi$}
    \popi^{xy}\leq 0 \Leftrightarrow \left\{\begin{aligned}
            &\langle T(x)-T(y),x-y\rangle \geq \mu \|x-y\|^2\\
            &\langle T(x)-T(y),x-y\rangle \geq \beta \|T(x)-T(y)\|^2,
    \end{aligned}\right.\label{eq:operators2}
\end{equation}
where again we use the notation $q$ instead of $p$ to denote a property.

Similarly to the case of monotone Lipschitz operators, \eqref{eq:operators2} is a pairwise interpolation condition on $\R^d$, and a global interpolation condition when $d=1$. However, when $K\geq3$ or $d\geq 2$, \eqref{eq:operators} is only necessary for interpolation and is no longer a $K$-wise interpolation condition. \cite[Proposition 4]{ryu2020operator}.
\subsection{Derivation of the one-point strengthening}
We thus strengthen \eqref{eq:operators2}, starting with a set of cardinality $3$.
\begin{proposition} \label{prop:1pt_strengthening_operators2}
    Consider $Q=\{(x_1,x_2,x_3)\}$, and \eqref{eq:operators2}.
    Then, $\tpopi$, the terwise one-point strengthening of \eqref{eq:operators2} with respect to $\R^d$, is satisfied only if
    \begin{equation}
        \tpopi^{ijk}\leq 0 \Rightarrow\left\{\begin{aligned}
                &\forall  \lambda_k, \ \mu_i, \ \mu_j \geq 0,\\
                0&\geq \lambda_k(\mu_i+\mu_j)(B_{ik}\mu_i+B_{jk}\mu_j)-\lambda_k \mu_i\mu_j B_{ij}+\lambda_k^2 (A_{ik} \mu_i+A_{jk}\mu_j)\\
                \MoveEqLeft[-1]+A_{ij} \mu_i\mu_j(\mu_i+\mu_j+2\lambda_k-4\beta\mu \lambda_k), \\
                0&\geq \lambda_k(\mu_i+\mu_j)(A_{ik}\mu_i+A_{jk}\mu_j)-\lambda_k \mu_i\mu_j A_{ij}+\lambda_k^2 (B_{ik} \mu_i+B_{jk}\mu_j)\\
                \MoveEqLeft[-1]+B_{ij} \mu_i\mu_j(\mu_i+\mu_j+2\lambda_k-4\beta\mu \lambda_k).
        \end{aligned}\right.\label{eq:final_popi}
    \end{equation}
    where $A_{ij}=-\langle t_i-t_j,x_i-x_j\rangle+\mu \|x_i-x_j\|^2$ and $B_{ij}=-\langle t_i-t_j,x_i-x_j\rangle+\beta \|t_i-t_j\|^2$.
\end{proposition}
\begin{proof}
    The terwise one-point strengthening of \eqref{eq:operators2} on $\R^d$ is given by
    \begin{align*}
        \tpopi= \max_{z\in \R^d}\min_{\substack{\tau \in \R, \\ t_z\in \R^d}} &\tau  && \\
        \text{s.t.} &-\langle t_z-t_i,z-x_i\rangle+\beta\|t_z-t_i\|^2\leq \tau,\quad \quad \quad \quad i\in \llbracket 3\rrbracket  && (\lambda_i) \nonumber \\
        & -\langle t_z-t_i,z-x_i\rangle+\mu\|z-x_i\|^2\leq \tau, \quad \quad \quad \quad  i\in \llbracket 3\rrbracket,&& (\mu_i) \nonumber
    \end{align*}
    where $\lambda_i,\mu_i$ are the associated dual coefficients. Partial dualization of the inner problem over $\tau$, and minimization over $t_z$ yields
    \begin{align*}
        \tpopi\geq \max_{\substack{z\in \R^d,\\ \lambda_i, \ \mu_i \geq 0}} \ &\sum_{i}\bigg  (\lambda_i\big(-\langle t_z-t_i,z-x_i\rangle+\beta\|t_z-t_i\|^2\big) +\mu_i \big(-\langle t_z-t_i,z-x_i\rangle+\mu\|z-x_i\|^2 \big)\bigg)\\
        \text{ s.t. }\ \  &\sum_i 2\beta \lambda_i (t_z^\star-t_i)-(\lambda_i+\mu_i) (z-x_i)=0, \sum_i (\lambda_i+\mu_i) =1,
    \end{align*}
    where $t_z^\star$ minimizes the inner problem in $\tpopi$. Minimizing over $z$, multiplying by $1-4\beta \mu(\sum_i \lambda_i)(\sum_i \mu_i)$ and rearranging the terms yields
    \begin{align*}
        0\geq \max_{ \substack{\lambda_i, \ \mu_i \geq 0\\4\beta \mu(\sum_i \lambda_i)(\sum_i \mu_i)\big)\leq 1 \\ \sum_i (\lambda_i+\mu_i) =1}}\ \sum_{i,j \in \llbracket 3 \rrbracket:\ j>i} \bigg(&B_{ij} \big(\lambda_i\lambda_j\sum_k \lambda_k+\lambda_i \mu_j (\sum_{k} \mu_k)+\lambda_j \mu_i (\sum_{k} \mu_{k})\\[-0.5cm]
            \MoveEqLeft[-4]+(2-4\beta \mu)(\sum_k\mu_k)\lambda_i\lambda_j-(\sum_{k} \lambda_{k})\mu_i\mu_j\big)\\\MoveEqLeft[-1] +A_{ij}\big(\mu_i\mu_j\sum_k \mu_k+(\mu_i \lambda_j (\sum_{k} \lambda_k)+\mu_j \lambda_i (\sum_{k} \lambda_{k})\\
        \MoveEqLeft[-4]+(2-4\beta \mu)(\sum_k\lambda_k)\mu_i\mu_j-(\sum_{k} \mu_{k})\lambda_i\lambda_j)\big)\bigg).
    \end{align*}
    This condition is homogeneous of degree $3$ in $\lambda_k, \ \mu_i,\ \mu_j$. Hence, for any $K\geq 0$, it must hold with $\tilde \mu_i=K\mu_i$, $\tilde \lambda_i=K\lambda_i$ replacing $\mu_i$, $\lambda_i$ respectively, which allows removing the conditions $\sum_i (\lambda_i+\mu_i) =1$ and $4\beta \mu(\sum_i \lambda_i)(\sum_i \mu_i)\big)\leq 1$.

    Setting only $\lambda_k,\ \mu_i$ and $\mu_j$ non-zero yields the first condition in \eqref{eq:final_popi}, and setting only $\mu_k,\ \lambda_i$ and $\lambda_j$ non-zero yields the second condition in \eqref{eq:final_popi}, hence $\tpopi$ is at least as strong as \eqref{eq:final_popi}.
\end{proof}
\begin{remark}
    Setting $\lambda_k=0$ in \eqref{eq:final_popi} allows recovering \eqref{eq:operators2}, hence again this relaxation of $\tpopi$ is at least as strong as \eqref{eq:operators2}. Given the example illustrated in \cref{fig:operators2}, \eqref{eq:final_popi} is actually strictly stronger than \eqref{eq:operators2}.
\end{remark}
\subsection{Derivation of a SOS reformulation}
\begin{figure}
    \centering
    \includegraphics[width=0.5\linewidth]{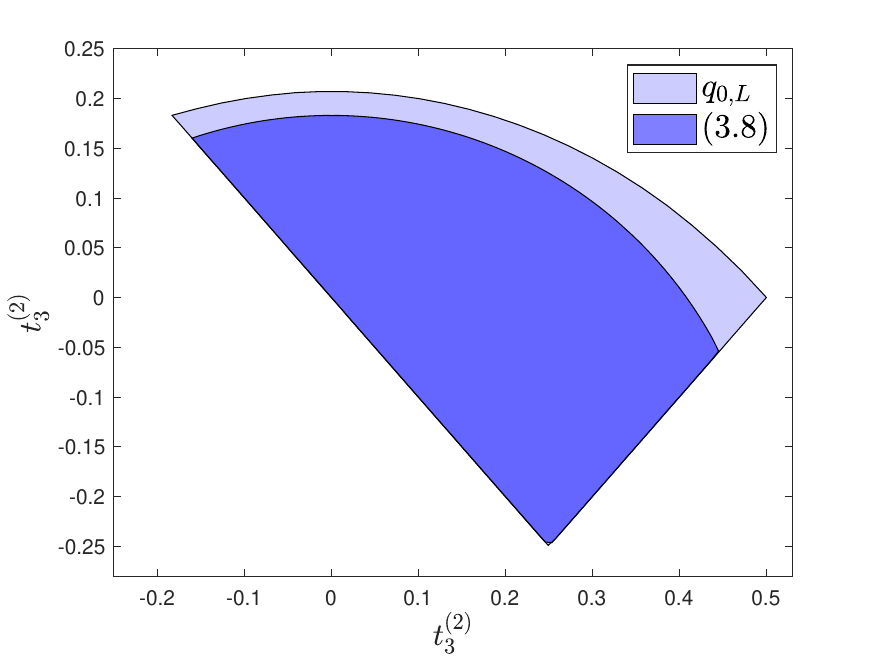}
    \caption{Allowed regions for $t_3$, according to \eqref{eq:operators2} or \eqref{eq:final_popi}, given $\mu=\frac{1}{2},\ \beta=1$, $(x_1,t_1)=\left(\binom{0}{0},\binom{0}{0}\right)$, $(x_2,t_2)=\left(\binom{1}{-1/2},\binom{1}{0}\right)$, $x_3=\binom{1/2}{1/2}$.}
    \label{fig:operators2}
\end{figure}

Similarly to \cref{sec:operators}, we propose a relaxation of \eqref{eq:final_popi} as an SDP.
\begin{proposition}\label{prop:SOS_opi}
    With the same notation as in \cref{prop:1pt_strengthening_operators2} and considering a triplet $S=\{(x_i,t_i)\}_{i\in \llbracket 3\rrbracket}$, it holds that $S$ satisfies \eqref{eq:final_popi} if and only if ~$\forall i,j,k \in \llbracket 3\rrbracket$, there exists some $M_{A,B}^{ijk}\succeq0$, $M_{B,A}^{ijk}\succeq0$, where
    \begin{equation}
        M_{A,B}=\begin{bmatrix}\label{def:matrix_opi}
            -B_{ij} & 0 & 0 & M_{14} &M_{15}& M_{16}&M_{17}\\
            0 & -A_{jk} & M_{23} & 0 & M_{25}&M_{26}&M_{27}\\
            0 & M_{23} & -B_{ij}  & M_{34}& M_{35}&0&M_{37}\\
            M_{14} & 0 & M_{34} &-B_{jk}&M_{45}&M_{46}&0\\
            M_{15}&M_{25} &M_{35} & M_{45}&M_{55}&M_{56}&M_{57}\\
            M_{16}&M_{26} &0 &M_{46} & M_{56}&-A_{ik}&0\\
            M_{17}&M_{27} &M_{37} &0 &M_{57} & 0& -B_{ik}
        \end{bmatrix},
    \end{equation}
    and
    \begin{align*}
        &M_{14}+M_{25}=
        M_{15}+M_{23}=
        M_{16}+M_{35}=
        M_{27}+M_{45}=
        M_{37}+M_{56}=
        M_{46}+M_{57}=0\\
        &M_{55}+2M_{17}+2M_{26}+2M_{34}=A_{ij}-A_{jk}-A_{ik}-(2-4\beta \mu)B_{ij}.
    \end{align*}
\end{proposition}
\begin{proof}
    The proof follows the same lines as that of \cref{prop:SOS_op}.
\end{proof}
We illustrate in \cref{fig:operators2} the difference in datasets satisfying  \eqref{eq:operators2} as compared to \eqref{eq:final_popi}.

\section{Uniformly convex functions}\label{sec:unif_conv}
This section studies the class $\F_{\mu,q}(\R^d)$ of uniformly convex functions $f$ satisfying, $\forall x,y \in \R^d$ \cite{iouditski2014primal},
\begin{equation}\tag{$p_{\mu,q}$}
    p_{\mu,q}^{xy}\leq 0\Leftrightarrow
    f(x)\geq f(y)+\langle \nabla f(y), x-y\rangle+\tfrac{\mu}{q}\|x-y\|^q, \text{where $\mu\geq 0, \ q\geq 2$}, \label{eq:unif_convex}
\end{equation}

We first prove that \eqref{eq:unif_convex} serves as a pairwise extension condition for $\F_{\mu,q}(\R^d)$ .
\begin{proposition}\label{prop:interp_unif_convex}
    Let $\F_{\mu,q}(\R^d)$. Then, \eqref{eq:unif_convex} is a pairwise extension condition for $\F_{\mu,q}(\R^d)$.
\end{proposition}
\begin{proof}
    Consider a set $S=\{(x_i,f_i,g_i)\}_{i =1,2}$ satisfying \eqref{eq:unif_convex}. We show that, given any arbitrary $x\in \R^d$, \eqref{eq:unif_convex} can be extended to $x$. To this end,  we propose quantities $f\in \R, g\in \R^d$ extending \eqref{eq:unif_convex} to $x$, that is such that $\forall i,j =1,2$, it holds
    \begin{align}
        f\geq &f_i+\langle g_i,x-x_i\rangle+\frac{\mu}{q}\|x-x_i\|^{q},\label{goal1}\\
        f_j\geq&f+\langle g,x_j-x\rangle+\frac{\mu}{q}\|x-x_j\|^{q}.\label{goal2}
    \end{align}
    Let
    \begin{align*}
        i_\star&=\underset{i=1,2}{\text{argmax}} \big(f_i+\langle g_i,x-x_i\rangle+\frac{\mu}{q}\|x-x_i\|^{q}\big), \\
        f&=f_{i_\star}+\langle g_{i_\star},x-x_{i_\star}\rangle+\frac{\mu}{q}\|x-x_{i_\star}\|^{q}, \
        g=g_{i_\star}+\mu\|x-x_{i_\star}\|^{q-2}(x-x_{i_\star})+K(x-x_{i_\star}),
    \end{align*}
    where $K\geq 0$.
    Then, \eqref{goal1} is automatically satisfied. In addition, it holds
    \begin{equation*}
        f+\langle g,x_{i_\star}-x\rangle+\frac{\mu}{q}\|x-x_{i_\star}\|^{q}=f_{i_\star}-\mu \|x-x_{i_\star}\|^{q}+ 2\frac{\mu}{q}\|x-x_{i_\star}\|^{q}-K\|x-x_{i_\star}\|\leq  f_{i_\star},
    \end{equation*}
    and
    \begin{equation*}
        \begin{aligned}
            f+\langle g,x_{j}-x\rangle+\frac{\mu}{q}\|x-x_{j}\|^{q}
            &=f_{i_\star}+\langle g_{i_\star},x_j-x_{i_\star}\rangle+ \langle K (x-x_{i_\star}),x_{j}-x\rangle\\
            \MoveEqLeft[-1]+\mu\|x-x_{i_\star}\|^{q-2}\langle x-x_{i_\star},x_j-x\rangle+\frac{\mu}{q}(\|x-x_{i_\star}\|^{q}+\|x-x_{j}\|^{q})\\
            &\leq f_j- \langle K X_{i_\star},X_j\rangle\\
            \MoveEqLeft[-1]+\frac{\mu}{q}(\| X_{i_\star}\|^{q}+\|X_j\|^{q}-\| X_{i_\star}-X_j\|^q-q\| X_{i_\star}\|^{q-2}\langle  X_{i_\star},X_j\rangle),
        \end{aligned}
    \end{equation*}
    where $ X_{i_\star}=x-x_{i_\star}$ and $X_j=x-x_j$. Suppose first $\langle  X_{i_\star},X_j\rangle\leq 0$, and let $K=0$. Then, by Bernoulli's inequality, and since $(x+y)^q\geq x^q+y^q$ for $q\geq 2, \ x,y\geq 0$, it holds
    \begin{equation*}
        \begin{aligned}
            \| X_{i_\star}-X_j\|^q&=(\| X_{i_\star}\|^2+\|X_j\|^2-2\langle X_{i_\star},X_j\rangle )^{q/2}\\
            &=(\| X_{i_\star}\|^2+\|X_j\|^2)^{q/2}(1-2\frac{\langle X_{i_\star},X_j\rangle}{\| X_{i_\star}\|^2+\|X_j\|^2} )^{q/2}\\
            &\geq\| X_{i_\star}\|^q+\|X_j\|^q-q\langle X_{i_\star},X_j\rangle(\| X_{i_\star}\|^2+\|X_j\|^2)^{q/2-1}\\
            &\geq \| X_{i_\star}\|^q+\|X_j\|^q-q\langle X_{i_\star},X_j\rangle \| X_{i_\star}\|^{q-2},
        \end{aligned}
    \end{equation*}
    hence \eqref{goal2} holds. Suppose now $\langle  X_{i_\star},X_j\rangle> 0$. Then there exists a sufficiently large $K$ such that \eqref{goal2} holds.
\end{proof}

Relying on the refinement procedure of \cref{sec:formal_extension}, we now show that \eqref{eq:unif_convex} can be strengthened by considering sets of higher cardinality.
\begin{proposition} \label{prop:1pt_strengthening_unif_convex}
    Consider $Q=\{(x_1,x_2,x_3)\}$, and \eqref{eq:unif_convex}. Then, $\tilde{p}_{\mu,q}$, the terwise one-point strengthening of \eqref{eq:unif_convex} on ~$\R^d$, is satisfied only if
    \begin{equation}
        \begin{aligned}[t]
            \tilde{p}_{\mu,q}^{ijk}\leq 0 \Rightarrow
            0\geq &\max_{\lambda \in [0,1]} \ \lambda (-f_i+f_k+\langle g_k,x_i-x_k\rangle)+(1-\lambda) (-f_j+f_k+\langle g_k,x_j-x_k\rangle)\\&+\frac{\mu}{q} \left((\lambda(1-\lambda)^q+(1-\lambda)\lambda^q)\|x_i-x_j\|^q+\|x_k-\lambda x_i-(1-\lambda)x_j\|^q\right).
        \end{aligned}\label{unif_conv_strengthened}
    \end{equation}
\end{proposition}

\begin{proof}
    The terwise one-point strengthening of \eqref{eq:unif_convex} on $\R^d$ is given by
    \begin{align*}
        \tilde{p}_{\mu,q}= \max_{z\in \R^d}\min_{\substack{\tau \in \R, \\ g_z\in \R^d,\ f_z\in \R}} \tau
        \text{ s.t. }
        & -f_z+ f_i+\langle g_i,z-x_i\rangle+ \frac{\mu}{q}\|z-x_i\|^q\leq \tau, \ i\in \llbracket 3\rrbracket && (\lambda_i)\nonumber\\
        &-f_i+ f_z-\langle g_z,z-x_i\rangle+  \frac{\mu}{q}\|z-x_i\|^q\leq \tau , \ i\in \llbracket 3\rrbracket&&(\mu_i)\nonumber
    \end{align*}
    where $\lambda_i, \mu_i$ ($i \in \llbracket 3 \rrbracket$) are the associated dual coefficients. Partial dualization of the inner problem over $\tau$ and $f_z$ yields
    \begin{align*}
        \tilde{p}_{\mu,q}= \max_{\substack{z\in \R^d,\\ \lambda_i, \ \mu_i\geq 0} }\ &\sum_{i=1}^3\big(\lambda_i (f_i+\langle g_i,z-x_i\rangle+ \frac{\mu}{q}\|z-x_i\|^q) + \mu_i (-f_i-\langle g_z,z-x_i\rangle+  \frac{\mu}{q}\|z-x_i\|^q)\big)\\
        \text{ s.t. }&
        \sum_{i=1}^3 \lambda_i= \sum_{i=1}^3 \mu_i=\frac{1}{2}, \quad z =2\sum_{i=1}^3(\mu_i x_i)\nonumber.
    \end{align*}
    In particular, setting only $\lambda_k=\frac{1}{2}$, $\mu_i, \mu_j=\frac{1}{2}-\mu_i$ non-zero yields a feasible solution of
    \begin{equation*}
        \begin{aligned}
            0&\geq \mu_i (-f_i+f_k+\langle g_k,x_i-x_k\rangle)+(\frac{1}{2}-\mu_i) (-f_j+f_k+\langle g_k,x_j-x_k\rangle)\\
            \MoveEqLeft[-1]+\frac{\mu}{q} \left((\mu_i(\frac{1}{2}-\mu_i)^q+(\frac{1}{2}-\mu_i)\mu_i^q)\|x_j-x_i\|^q+\|x_k-\mu_i x_i-(\frac{1}{2}-\mu_i)x_j\|^q\right),  \forall \mu_i\in [0,\frac{1}{2}].
        \end{aligned}
    \end{equation*}
    Hence $\tilde{p}_{\mu,L}$ is at least as strong as \eqref{unif_conv_strengthened}.
\end{proof}

\begin{remark}
    Setting $\lambda=1$ in \eqref{unif_conv_strengthened} allows recovering \eqref{eq:unif_convex}.
\end{remark}
\begin{remark}
    The case $p=2$ can be simplified as
    \begin{align*}
        \tilde{p}_{\mu,2}^{ijk}\leq 0 &\Rightarrow
        0\geq \lambda (-f_i+f_k+\langle g_k,x_i-x_k\rangle+\frac{\mu}{2}\|x_i-x_k\|^2)\\
        \MoveEqLeft[-3] +(1-\lambda) (-f_j+f_k+\langle g_k,x_j-x_k\rangle+\frac{\mu}{2}\|x_i-x_k\|^2)\\
        &\Leftrightarrow  f_i\geq f_k+\langle g_k,x_i-x_k\rangle+\frac{\mu}{2}\|x_i-x_k\|^2 \text{ and } f_j\geq f_k+\langle g_k,x_j-x_k\rangle+\frac{\mu}{2}\|x_j-x_k\|^2.
    \end{align*}
    We recover the classical definition of $\mu$-strongly convex functions.
\end{remark}

\section{Smooth convex functions with constrained convex domain}\label{sec:extension_subset}

Given an open convex subset $\bar Q\subsetneq \R^d$, consider the class $\F_{0,L}(\bar Q)$ of functions that are smooth and convex on $\bar Q$, in the sense of \eqref{eq:convexity_gradient_lipschitz}, but not necessarily outside of $\bar Q$. It is known that \eqref{eq:true_smooth_convexity} no longer holds in general for this function class \cite[Section 2]{drori2018properties}, but that it remains valid for pairs $(x,y)$ of points $x,y \in \bar Q$ that are sufficiently close. In particular, it holds at every $x,y\in \bar Q$ such that \cite[Theorem 3.1]{drori2018properties} \begin{equation}
    \|x-y\|< \mathrm{dist}(y,\R^d\setminus \bar Q) ,\label{eq:cond_drori}
\end{equation}

This result is obtained by showing that in this case, all quantities involved in the strengthening of \eqref{eq:convexity_quadratic_UB}  into \eqref{eq:true_smooth_convexity} belong to $\bar Q$. Building on the iterative strengthening of \eqref{eq:convexity_gradient_lipschitz}, we improve this description by showing \eqref{eq:true_smooth_convexity} holds at pairs that are twice as distant.
\begin{proposition}\label{prop:cond_drori_improved}
    Let $\bar Q\subsetneq \R^d$an open convex set and $f\in \F_{0,L}(\bar Q)$. Then, for any $x,y\in \bar Q$ such that \begin{equation}
        \|x-y\|< 2\mathrm{dist}(y,\R^d\setminus \bar Q) ,\label{eq:cond_drori_improved}
    \end{equation}
    it holds
    \begin{equation}
        f(y)\geq f(x)+\langle \nabla f(x),y-x\rangle +\tfrac{1}{2L}\|\nabla f(x)-\nabla f(y)\|^2.
    \end{equation}
\end{proposition}
\begin{proof}
    Let $z=y+\frac{g(x)-g(y)}{2L}$. It holds
    \begin{equation*}
        \|z-y\|=\tfrac{1}{2L}\|g(y)-g(x)\|\leq \tfrac{1}{2}\|y-x\|\underset{\eqref{eq:cond_drori_improved}}{<} \mathrm{dist}(y,\R^d\setminus \bar Q),
    \end{equation*}
    hence $z\in \bar Q$. Setting $x_2=x$ and $x_1=y$ in the proof of \cref{prop:1pt_strengthening_smooth_convex}, observe that the strengthening of \eqref{eq:convexity_gradient_lipschitz} into \eqref{eq:true_smooth_convexity} remains valid since it only involves quantities belonging to $\bar Q$.
\end{proof}

Building on \cite[Theorem 3.1]{drori2018properties}, the author of \cite{drori2018properties} obtained properties valid for all $x,y\in \bar Q$, however far from each other, that consists of bounds on $f(y)-f(x)$, solutions to given optimization problems. Relying on \cref{prop:cond_drori_improved} instead of \cite[Theorem 3.1]{drori2018properties} allows improving this description of $\F_{0,L}(\bar Q)$.
\begin{corollary}[from \protect{\cite[Corollary 3.2]{drori2018properties}}]\label{cor:drori_improved}
    Given $\{(x,g_x,f_x),(y,g_y,f_y)\}\in (\R^d\times \R^d \times \R)^2$, let \begin{align*}
        \mathcal{U}_N (\text{resp. }\mathcal{B}_N)(x,g_x,f_x,y,g_y)=&\max (\text{resp. }\min)_{g_i\in \R^d, \ f_i\in \R}\quad  f_N&&\\
        \text{s.t. } &\frac{1}{2L} \|g_i-g_{i+1}\|^2\leq f_i-f_{i+1}-\frac{1}{N} \langle g_{i+1},x-y\rangle, &&0\leq i<N\\
        &\frac{1}{2L} \|g_i-g_{i+1}\|^2\leq f_{i+1}-f_{i}-\frac{1}{N} \langle g_{i},y-x\rangle, &&0\leq i<N\\
        &f_0=f_x, \ g_0=g_x, \ g_N=g_y.&&
    \end{align*}
    Then,
    \begin{enumerate}[label=(\roman*)]
        \item If $f\in \F_{0,L}(\bar Q)$, then for any $x,y\in \bar Q$, and $\forall N>\frac{\|x-y\|}{2\min(\mathrm{dist}(x,\R^d\setminus \bar Q),\mathrm{dist}(y,\R^d\setminus \bar Q))}$,
            \begin{equation*}
                \mathcal{B}_N(x,\nabla f(x),f(x), y,\nabla f(y))\leq f(y)\leq \mathcal{U}_N(x,\nabla f(x),f(x), y,\nabla f(y)).
            \end{equation*}
        \item If $\{(x,g_x,f_x),(y,g_y,f_y)\}\in (\R^d\times \R^d \times \R)^2$ satisfy, for some $L,N>0$,
            \begin{equation*}
                \mathcal{B}_N(x,g_x,f_x,y,g_y)\leq f_y\leq \mathcal{U}_N(x,g_x,f_x,y,g_y),
            \end{equation*}
            then there exists a function $f\in \F_{0,L}(\bar Q)$ such that $f(x)=f_x,\ f(y)=f_y, \ \nabla f(x)=g_x, \text{ and } \nabla f(y)= g_y$.
    \end{enumerate}
\end{corollary}
\begin{proof}
    The proof follows the same lines as in \cite[Corollary 3.2]{drori2018properties}, replacing Proposition \cite[Theorem 3.1]{drori2018properties} by \cref{prop:cond_drori_improved}.
\end{proof}
In \cite[Corollary 3.2]{drori2018properties}, a similar result is proved, except for the factor $2$ in the definition of $N$: \cref{cor:drori_improved} allows bounding $f(y)$ by $\mathcal{U}_N$, $\mathcal{B}_N$ for twice smaller values of $N$ than  \cite[Corollary 3.2]{drori2018properties}.

\begin{figure}[t]
    \centering
    \includegraphics[width=0.5\linewidth]{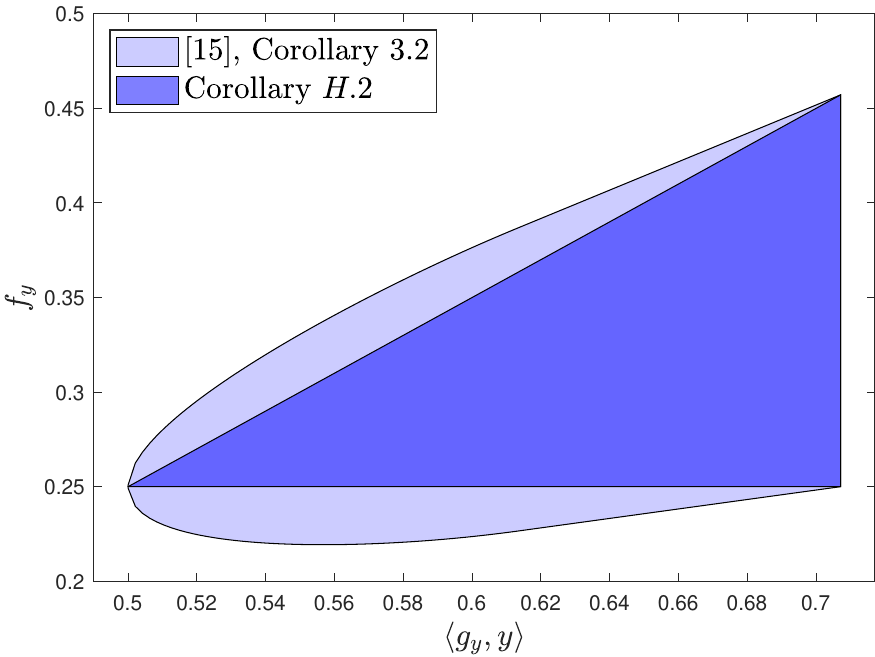}
    \caption{From \cite[Figure 3]{drori2018properties}. Allowed regions for $f(y)$ as a function of $\langle \nabla f(y),y\rangle $, according to \cite[Corollary 3.2]{drori2018properties} and \cref{cor:drori_improved}, and given $L=1$, $(x,f(x),\nabla f(x))=(0,0,0), \ \|y\|^2=1$, $\|\nabla f(y)\|^2=\frac{1}{2}$ and $\min(\mathrm{dist}(x,\R^d\setminus \bar Q),\mathrm{dist}(y,\R^d\setminus \bar Q))=1$. \Cref{cor:drori_improved} allows bounding $f_y$ by $\mathcal{B}_1$, $\mathcal{U}_1$, while \cite[Corollary 3.2]{drori2018properties} can only rely on $\mathcal{B}_2$, $\mathcal{U}_2$.}
    \label{fig:drori_improved}
\end{figure}

To illustrate how this strengthens the description of $\F_{0,L}(\bar Q)$, consider the following numerical example, from \cite[Figure 3]{drori2018properties}. Given $f:\bar Q\to \R$, suppose $L=1$, $x=0, \ \nabla f(x)=0, \ f(x)=0, \ \|y\|^2=1$, $\|\nabla f(y)\|^2=\frac{1}{2}$ and $\min(\mathrm{dist}(x,\R^d\setminus \bar Q),\mathrm{dist}(y,\R^d\setminus \bar Q))=1$. We show in \cref{fig:drori_improved} the allowed regions for $f(y)$ as a function of $\langle \nabla f(y),y\rangle $, according to \cite[Corollary 3.2]{drori2018properties} and \cref{cor:drori_improved}.

\bibliographystyle{jnsao}
\bibliography{main.bib}
\end{document}